\documentclass[10pt]{article} 

\usepackage[colorlinks=true, allcolors=blue]{hyperref}

\usepackage[utf8]{inputenc} 

\usepackage{geometry} 
\usepackage{graphicx} 

\usepackage[round]{natbib}

\usepackage{booktabs} 
\usepackage{array} 
\usepackage{paralist} 
\usepackage{verbatim} 
\usepackage{subfig} 

\usepackage{amsmath}
\usepackage{amsfonts}
\usepackage{amssymb}
\usepackage{dsfont}
\usepackage{amsthm}

\usepackage{multirow}

\allowdisplaybreaks

\theoremstyle{theorem}
\newtheorem{prop}{Proposition}
\newtheorem{thm}{Theorem}
\newtheorem{lem}{Lemma}
\newtheorem{cor}{Corollary}
\theoremstyle{definition}
\newtheorem{defn}{Definition}

\newcommand{\RRR}{\mathbb{R}}
\newcommand{\EEE}{\mathbb{E}}
\newcommand{\PPP}{\mathbb{P}}
\newcommand{\dd}{\mathrm{d}}
\newcommand{\one}{\mathds{1}}
\newcommand{\boldw}{\mathbf{w}}

\newcommand{\boldP}{\mathbf{P}}

\newcommand{\ds}{\displaystyle}

\newcommand{\GDS}{\mathrm{GDS}}
\newcommand{\GDSd}{\mathrm{GDS}_\mathrm{d}}
\newcommand{\GDSc}{\mathrm{GDS}_\mathrm{c}}
\newcommand{\LEPSr}{\mathrm{LEPS}_\mathrm{r}}
\newcommand{\LEPSc}{\mathrm{LEPS}_\mathrm{c}}
\newcommand{\Grin}{\mathrm{GRIN}}
\newcommand{\SEEPS}{\mathrm{SEEPS}}

\newcommand{\cN}{\mathcal{N}}
\newcommand{\RG}{\mathcal{R}(G)}

\newcommand{\twocdots}{\,\cdot\,,\cdot\,}

\usepackage{fancyhdr} 
\usepackage{sectsty}
\allsectionsfont{\sffamily\mdseries\upshape} 

\usepackage[nottoc,notlof,notlot]{tocbibind} 
\usepackage[titles,subfigure]{tocloft} 

\title{On equitable scoring functions and optimal forecasting behaviour}

\author{Robert J. Taggart\thanks{Corresponding author: robert.taggart@bom.gov.au} \quad and \quad Nicholas Loveday\\[1ex] Australian Bureau of Meteorology}

\date{6 August 2026}

\begin{document}

\maketitle


\begin{abstract}
\noindent Equitable scoring functions have a long history in meteorological forecast verification and have recently gained renewed prominence through the use of the Stable Equitable Error in Probability Space (SEEPS) score for evaluating precipitation forecasts from numerical and machine-learning weather prediction systems.

This paper provides a systematic analysis of equitable scoring functions through the optimal forecasts they induce. We introduce the \textit{generalized diagonal score}, which defines a broad class of equitable scoring functions that includes, up to equivalence, the SEEPS, Gerrity, Peirce, Barnston and diagonal scores. Within this framework, we show that optimal single-valued and categorical forecasts are characterized by \textit{crossing points} between the predictive and climatological cumulative distribution functions. Moreover, we show that the generalized diagonal score, when evaluating predictive distributions, is proper but not strictly proper, and is insensitive to substantial forecast misspecification away from crossing points. The framework also yields scoring rules that are both proper and equitable for probabilistic forecasts with categorical outcomes.

For single-valued and categorical forecasts, the crossing-point characterization shows that optimal forecasts can vary across climatologies even when the underlying predictive distribution is unchanged. Consequently, identical predictive distributions may lead to substantially different optimal forecasts under equitable scores, with important implications for assessing the suitability of such scores for any given application.
\end{abstract}

\vspace{10pt}


\section{Introduction}\label{s:intro}

Forecast evaluation plays a central role in ranking forecasting systems and improving forecast processes. A common approach is to use a scoring function, which assigns a numerical value (the score) to a forecast given the corresponding observation. When the score is negatively oriented, it represents a penalty, and the goal is to minimize this value. The choice of scoring function is therefore critical, as it determines not only how forecast performance is assessed but also the type of forecasting behaviour that is encouraged.

When forecasts are probabilistic, scoring functions are typically referred to as scoring rules. A broad multidisciplinary consensus holds that scoring rules should be \textit{proper}, meaning that a forecaster minimizes their expected score by issuing a predictive distribution that reflects their true beliefs \citep{murphy1967note, gneiting2007strictly}. Proper scoring rules therefore promote honest forecasting and ensure that available information is used effectively.

In many applications, however, forecasts are not probabilistic but instead take the form of single-valued\footnote{We use the term \textit{single-valued forecast} for a forecast that consists of a single real number. Such forecasts are sometimes called deterministic forecasts in meteorology and point forecasts in the statistical forecasting literature. The term point forecast is avoided here because it may be confused with forecasts at an individual spatial location in contrast to those on a grid.} or categorical predictions. In such settings, it is natural to specify a \textit{directive} describing how a forecaster should convert probabilistic beliefs into a forecast, such as issuing the expected value of the predictive distribution or forecasting an event when its probability exceeds a given threshold. A scoring function is said to be \textit{consistent} with a directive if the forecaster's expected score is minimized by following it \citep{murphy1985forecast, gneiting2011making}. When consistency holds, there is no conflict between optimizing the score and adhering to the directive. The forecast that minimizes the forecaster's expected score is called the Bayes act \citep{gneiting2011making}. Bayes acts provide a useful way to characterize the forecasting behaviour induced by a scoring function.

Alongside propriety and consistency, another property known as \textit{equitability} has been widely used in the meteorological literature. A scoring function is equitable if all constant forecasts receive the same expected score under the climatological distribution. The term `equitable' appears to have been coined by \citet{gandin1992equitable}, though the core idea is evident in the work of \citet{gringorten1951verification, gringorten1965measure, gringorten1967verification} and is a property possessed by the \citet{peirce1884numerical} skill score. Though several authors have argued that equitability is a desirable property for a variety of reasons \citep{gandin1992equitable, potts1996revised, Livezey2003Categorical, jolliffe2008impenetrable, hogan2010equitability, rodwell2010new}, its implications for forecast behaviour have not been fully appreciated. In particular, with a few exceptions (e.g., \citealt[Fig.~3]{gandin1992equitable}; \citealt{ben2021verification}), it has generally not been made explicit what types of forecasts are encouraged when a forecaster seeks to optimize their expected score under an equitable scoring function.

These questions are of increasing practical importance due to the growing use of equitable scores in operational and ML forecast evaluation. In particular, the Stable Equitable Error in Probability Space (SEEPS) score \citep{rodwell2010new} has been widely adopted for verifying precipitation forecasts from weather prediction systems, including both physics-based and data-driven models \citep{north2013assessment, lavers2022evaluation, lang2024aifs, rasp2024weatherbench, radford2025comparison, pan2025comparative, sha2025investigating, moldovan2025update}. Despite its widespread use, comparatively little is understood about the forecasting behaviour induced by optimizing SEEPS, or how this behaviour relates to the forecaster's underlying probabilistic judgments.

The purpose of this paper is to provide a systematic and mathematically rigorous analysis of equitable scoring functions, with a focus on the forecasting strategies that they induce. In particular, we characterize the Bayes acts associated with a broad class of equitable scores and examine their relationship to the underlying climatological distribution. This perspective provides new insight into how equitable scores interact with probabilistic beliefs, and into the types of probabilistic, single-valued and categorical forecasts that they encourage in practice.

\begin{table}[t!]
    \renewcommand{\arraystretch}{1.5} 
    \centering \small
    \begin{tabular}{|l|l|}
    \hline
	forecast type & score \\
   \hline
	probability distribution on $\RRR$ & diagonal score \citep{ben2018diagonal} \\
		& diagonal elementary score \citep{ben2018diagonal} \\
   \hline
	single-valued in $\RRR$ & score of \citet{gringorten1965measure} \\
		& revised LEPS score \citep{potts1996revised} \\
		& SEEPS \citep{rodwell2010new} (via discretization) \\
		& score of \cite{ben2021verification} \\
   \hline
	categorical & skill score of \citet{peirce1884numerical} \\
		& \citet{gerrity1992note} score \\
		& scores of \citet{barnston1992correspondence} \\
		& SEEPS \citep{rodwell2010new} \\
    \hline
    \end{tabular}
    \caption{Published equitable scores that are, up to equivalence, special cases of the generalized diagonal score.}
    \label{tab:scores in GDS form}
\end{table}

The main theoretical contributions of this paper (Sections~\ref{s:equitability}–\ref{s:optimal}) are as follows.
\begin{itemize}
\item We introduce, building on the work of \citet{ben2018diagonal}, the \emph{generalized diagonal score} (GDS), which defines a broad class of equitable scoring functions, including (up to equivalence) the scores listed in Table~\ref{tab:scores in GDS form}.
\item We give precise statements on the regularity properties of the climatological distribution such that equitability holds with GDS and the scores of Table~\ref{tab:scores in GDS form}.
\item We construct new equitable scores from GDS when the climatological distribution is an empirical cumulative distribution function (CDF).
\item We characterize the local minimizers of the expected generalized diagonal score for single-valued and categorical forecasts in terms of \emph{crossing points} (c.f. \citealt{ben2021verification}). Crossing points may be viewed informally as points where the difference between the forecaster's CDF and the climatological CDF changes sign from negative to positive. See Definition~\ref{def:crossing point} for a precise definition. This characterization implies that, when the crossing point is unique, it is an optimal forecast for all equitable scores in Table~\ref{tab:scores in GDS form} in the single-valued and categorical settings. When multiple crossing points exist, an optimal forecast lies among them, though which is optimal depends on parameters that specify the score.
\item We establish that the GDS, when scoring predictive distributions, is proper but not strictly proper. In particular, substantially misspecified forecasts can attain the same expected score as the forecaster's true predictive distribution whenever their quantiles lie on the same side of the corresponding climatological quantiles.
\item We construct, using GDS, equitable proper scoring rules for probabilistic forecasts with categorical outcomes, including the binary case, thereby showing that the apparent incompatibility between propriety and equitability disappears when the smoothness assumptions of \citet{jolliffe2008proper} are relaxed.
\item We derive explicit formulae for the LEPSCAT scoring matrix \citep{potts1996revised} in the $N$-category case with general climatological probabilities, and prove that LEPSCAT is equitable. We further show that, in general, LEPSCAT cannot be represented in GDS form and that its Bayes act cannot be characterized via crossing points. Consequently, although the Bayes act for many equitable scores admit a crossing point characterization, equitability alone does not guarantee such a characterization.
\end{itemize}

Although our analysis is framed in terms of predictive distributions and decision theory, it is also relevant in settings where forecasting systems are developed, selected, or tuned by optimizing verification scores. Because scoring functions induce particular optimal forecast behaviours, improvements in score can arise not only from improved predictive information but also from systematic shifts towards behaviours rewarded by the score.

The final two sections of the paper focus on the implications of this theory for our understanding of equitable scores, and in particular for the use of SEEPS in evaluating precipitation forecasts from global weather prediction models. Section~\ref{s:experiments} presents two experiments. The first illustrates the types of forecasting behaviour induced by equitable scores for precipitation forecasts. In particular, it shows how crossing-point optimality can cause SEEPS to favour substantially wetter forecasts in dry climates and substantially drier forecasts in wet climates, even when the underlying predictive distributions are identical. The second examines and challenges claims that equitable scores promote `unbiased' forecast behaviour, in the sense of aligning the marginal distributions of forecasts and observations. Section~\ref{s:discussion} then evaluates further claims about equitable scores in the literature, considers alternatives to SEEPS, and concludes with recommendations. Proofs of the new mathematical results are given in the appendix.

\section{Mathematical setting and notation}

Let $Y$ be the future quantity to be observed, such as daily maximum temperature or 24-hour precipitation accumulation at a specific location, which takes possible values from some interval $I$ contained in the real line $\RRR$. The cases of primary interest in meteorology are the real line $I=\RRR$, the nonnegative half-line $I=[0,\infty)$, and the positive half-line $I=(0,\infty)$. The variable $Y$ is regarded as a random variable and the corresponding observation, or realization of $Y$, is a specific value in $I$ denoted by $y$.

A predictive distribution for $Y$ is a probability distribution on $I$. To make precise statements in the paper, we adopt a measure-theoretic framework for probability distributions, but to make the body of the paper as accessible as possible, technical terms related to measure theory will be restricted to this paragraph, footnotes and the appendix. Suffice to say, the setting is very general. Let $\mathcal{P}$ denote the class of Borel probability measures with mass concentrated in $I$. Every predictive distribution $F$ for $Y$ is assumed to be a member of $\mathcal{P}$ and will be identified with its cumulative distribution function (CDF), such that $F(z)$ denotes the probability that $Y\leq z$ whenever $z\in \RRR$, and $F$ is right-continuous. We denote the left-hand limit of $F$ at $z$ by $F(z^-)$, so that $F(z^-)$ is the probability that $Y<z$. The class $\mathcal{P}$ includes distributions $F$ that have discontinuities, and hence accommodates variables such as precipitation, streamflow, and wind magnitude (which are typically mixtures of discrete and continuous components), as well as purely discrete distributions (such as empirical CDFs and distributions on categories). Throughout, $\mu$ denotes a Borel measure on $(0,1)$, not necessarily a probability measure. Whenever $a\in\RRR$, let $\delta_a$ denote the Dirac measure with unit mass concentrated at $a$.

Given $G$ in $\mathcal{P}$, define the restricted domain $\mathcal{D}(G)$ of $G$ by $\mathcal{D}(G)=\{z\in \RRR:0<G(z)<1\}$ and the restricted range $\mathcal{R}(G)$ of $G$ by $\mathcal{R}(G)=\{\alpha\in (0,1):G(z)=\alpha\text{ for some } z\in\RRR\}$.

If $F\in\mathcal{P}$, its generalized inverse (i.e., quantile function) $F^{-1}:(0,1)\to \RRR$ is defined by $F^{-1}(\alpha)=\min\{z\in \RRR: F(z)\geq \alpha\}$. Hence $F^{-1}(\alpha)$ is single-valued, even if there are multiple $z$ in $\RRR$ such that $F(z)=\alpha$. Note that, in general, $F^{-1}(F(z)) \leq z$ and $F(F^{-1}(\alpha)) \geq \alpha$ but equality is not guaranteed.

If $h:I\to\RRR$ is a function on $I$ and $F$ a probability distribution on $I$, then the expected value of $h(Y)$, given that $Y$ has distribution $F$, will be denoted $\EEE_F h(Y)$ when the expectation exists. Similarly, let $\PPP_F(Y>a)$ denote the probability that $Y>a$ given that $Y$ has distribution $F$.

The normal distribution with mean $m$ and variance $\sigma^2$ will be denoted by $\cN(m, \sigma^2)$.

The indicator function will be denoted by $\one$, so that, for example, $\one\{x>y\}$ equals $1$ if $x>y$ and equals $0$ otherwise.

\section{Equitability and equivalence}\label{s:equitability}

\citet{gandin1992equitable} state that ``equitable skill scores formulated in this paper are based on the fundamental concept that constant forecasts of any event, as well as forecasts produced by a procedure in which a forecast event is chosen at random, should be accorded the same level of skill.'' They then formalize this concept, in the context of categorical forecasting, by requiring that there exists a constant $c$ such that (i) the expected score for any forecast category, when observations follow the climatological distribution, equals $c$, and (ii) the expected score under a random forecasting strategy, when forecasts are uniformly distributed over the finite set of categories and observations follow the climatological distribution, also equals $c$ \citep[Eqs.~(2,3)]{gandin1992equitable}. This formulation is widely regarded as the starting point for definitions of equitability in the meteorological literature; see, for example, \citet{hogan2010equitability}.

We therefore adopt the following definition, which applies more generally than the categorical setting.

\begin{defn}\label{def:equitability}
Suppose that $S:\mathcal{F}\times\mathcal{O}\to [-\infty,\infty]$ is a scoring function on the prediction--observation domain $\mathcal{F}\times\mathcal{O}$, where forecasts take values in $\mathcal{F}$ and observations take values in $\mathcal{O}$. Let $G$ denote a probability distribution on $\mathcal{O}$.\footnote{The observation space $\mathcal{O}$ is assumed to be equipped with a suitable $\sigma$-algebra. In the applications considered here, $\mathcal{O}$ is either an interval $I$ with its Borel $\sigma$-algebra, or a finite set $\{C_1,\ldots,C_n\}$ of categories equipped with the power set $\sigma$-algebra.} We say that $S$ is \textit{equitable} with respect to $G$ if there exists a constant $c$ in $\RRR$ such that $\EEE_G S(x,Y) = c$ for every $x$ in $\mathcal{F}$. We call $c$ the \textit{equitability constant}.
\end{defn}

In meteorological applications, the distribution $G$ will usually be the climatological (i.e., marginal) distribution of $Y$. We make the dependence on $G$  in Definition~\ref{def:equitability} explicit because expectations can also be taken with respect to other distributions. Subject to suitable regularity conditions, if $S$ is equitable with respect to $G$, then a random forecasting strategy, where the forecast $X$ is drawn from some distribution $H$ on $\mathcal{F}$, also yields the same expected score $c$. Formally:
\begin{align*}
\EEE[S(X,Y) \mid X \sim H, Y \sim G]
&= \EEE[\EEE_G S(X,Y) \mid X\sim H] \\
&= \EEE[c \mid X \sim H] \\
&= c.
\end{align*}

In this paper, we primarily work with negatively oriented scoring functions $S:\mathcal{F}\times\mathcal{O}\to[0,\infty]$, where smaller values indicate better performance, with the additional location assumption that $S(x,y)\geq0$ for all $(x,y)$ in $\mathcal{F}\times\mathcal{O}$. However, some equitable scores in the literature, such as those of \citet{gandin1992equitable} and \citet{gerrity1992note}, are positively oriented, so that larger values indicate better performance. Scoring functions can be transformed to change orientation, scale, or location without altering their essential properties. The following definition follows \citet{gneiting2007strictly}.

\begin{defn}
Suppose that $S$ and $S'$ are two scoring functions on a prediction--observation domain $\mathcal{F}\times\mathcal{O}$. The scoring function $S'$ is said to be \textit{equivalent} to $S$ if there exist a nonzero constant $a$ and a function $h:\mathcal{O}\to\RRR$ such that
\begin{equation}
S'(x,y)=aS(x,y)+h(y) \qquad \forall x\in \mathcal{F},\,\,\forall y\in \mathcal{O}.
\end{equation}
\end{defn}

The constant $a$ is positive if $S$ and $S'$ have the same orientation, and negative otherwise. Under the assumption that the relevant expectations exist, it is straightforward to show that if $S$ and $S'$ are equivalent then
\begin{enumerate}
\item[(i)] $S$ is equitable if and only if $S'$ is equitable;
\item[(ii)] a forecast $x$ in $\mathcal{F}$ optimizes the expected score $\EEE_FS(x,Y)$ if and only if $x$ optimizes the expected score $\EEE_FS'(x,Y)$, where $F$ is a probability distribution on $\mathcal{O}$.
\end{enumerate}
Here, `optimizes' means `minimizes' if the score is negatively oriented, and `maximizes' if it is positively oriented.

\section{The generalized diagonal score}\label{s:GDS}
We introduce a generalized version of the diagonal score of \citet{ben2018diagonal}. This generalization is motivated by the need to develop a framework that encompasses a broad class of equitable scores, and to accommodate distributions that are not necessarily continuous. The construction builds on elementary scoring functions, which we briefly recall.

To begin, for $\alpha$ in $(0,1)$ and $\theta$ in $I$, let $S_{\alpha,\theta}:I\times I\to[0,1)$ denote the elementary scoring function of \citet{ehm2016quantiles}, defined by
\begin{align}
S_{\alpha,\theta}(x,y)
&=
\begin{cases}
\alpha & \text{if } x \leq\theta<y \\
1-\alpha & \text{if } x >\theta\geq y \\
0 & \text{otherwise}
\end{cases} \notag \\ 
&= (\one\{y<x\} - \alpha)(\one\{\theta < x\} - \one\{\theta < y\}),  \label{eq:elementary score}
\end{align}
where $x$ is the real-valued forecast in $I$ and $y$ is the real-valued observation in $I$. Here we follow the notational convention of \citet{ehm2016quantiles}, rather than that of \citet{ben2018diagonal} whose notation is recoverable from our presentation by replacing $\alpha$ with $1-\alpha$ in Equation~(\ref{eq:elementary score}). The score $S_{\alpha,\theta}(x,y)$ admits an interpretation within a binary decision model, where an event is the set of outcomes $y$ such that $y>\theta$ and a user takes preventative action if and only if an event is forecast, namely $x$ exceeds the decision threshold $\theta$. The score $S_{\alpha,\theta}(x,y)$ then assigns a penalty of $\alpha$ to a miss (event observed but not forecast), a penalty of $1-\alpha$ to a false alarm (event forecast but not observed), and a penalty of $0$ to correctly forecast events and non-events. As discussed by \citet{ehm2016quantiles} and \citet{ben2018diagonal}, the score $S_{\alpha,\theta}(x,y)$ is consistent with the classical cost–loss decision model \citep{murphy1977value, richardson2000skill}. An optimal forecast $x$, namely one that minimizes the forecaster's expected score, is an $\alpha$-quantile of the predictive distribution.

The generalized diagonal score is constructed from elementary scoring functions that are equitable. It is therefore natural to characterize the class of equitable elementary scoring functions.

\begin{prop}\label{prop:elementary score equitable}
Suppose that $G\in\mathcal{P}$, $\alpha\in(0,1)$ and $\theta\in I$. The elementary scoring function $S_{\alpha,\theta}:I\times I\to[0,1)$ is equitable with respect to $G$ if and only if $\alpha=G(\theta)$. The equitability constant is $\alpha(1-\alpha)$.
\end{prop}

In other words, $S_{\alpha,\theta}$ is equitable precisely when $\alpha$ equals the climatological probability of not exceeding the threshold $\theta$.

\begin{cor}\label{cor:diagonal score equitable}
Suppose that $G\in\mathcal{P}$ and $\alpha\in(0,1)$. Then  $S_{\alpha,G^{-1}(\alpha)}:I\times I\to[0,1)$ is equitable with respect to $G$ if and only if $\alpha\in\RG$. The equitability constant is $\alpha(1-\alpha)$.
\end{cor}

The scoring function $S_{\alpha,G^{-1}(\alpha)}$ is called the \textit{diagonal elementary score} and was noted by \citet{ben2018diagonal} to be equitable whenever $0<\alpha<1$ under the implicit assumption that $G$ is continuous. The restriction $\alpha \in \RG$ in Corollary~\ref{cor:diagonal score equitable} is essential when $G$ has discontinuities, which arises in important applications such as climatological distributions of wind magnitude and daily precipitation.

Equitable scoring rules can be constructed from diagonal elementary scores by averaging them over different values of $\alpha$. Given a predictive distribution $F$, an observation $y$ and a value $\alpha$ in $\RG$, one may score $F$ using the equitable diagonal elementary score $S_{\alpha,G^{-1}(\alpha)}(x,y)$, where the single-valued forecast $x$ is obtained from $F$ via $x=T_\alpha(F)$. Here $T_\alpha$ is a functional that maps distributions to elements of $I$. For example, taking $T_\alpha(F)=F^{-1}(\alpha)$ means that $T_\alpha$ is the $\alpha$-quantile functional. Averaging these elementary scores with respect to a measure $\mu$ on $(0,1)$ yields the scoring rule
\begin{equation}\label{eq:general equitable scoring rule}
S(F,y; T, G, \mu) = \int_{\RG} S_{\alpha,G^{-1}(\alpha)}(T_\alpha(F), y) \,\dd\mu(\alpha).
\end{equation}
Under suitable regularity conditions on the family $T$ of functionals, the measure $\mu$ and the class of distributions, $S$ is well-defined and is equitable with respect to $G$ by Corollary~\ref{cor:diagonal score equitable}.

Among the equitable scoring rules of the form \eqref{eq:general equitable scoring rule}, the quantile-based construction obtained by taking $T_\alpha(F)=F^{-1}(\alpha)$ plays a special role because it is both equitable and proper. Propriety follows from the consistency of $S_{\alpha,G^{-1}(\alpha)}$ for the $\alpha$-quantile functional \citep{ehm2016quantiles}, as recognised by \citet{ben2018diagonal} and developed further in Section~\ref{ss:optimal GDS}. We therefore focus on this case and refer to the resulting scoring rule as the generalized diagonal score.

The \textit{generalized diagonal score} $\GDS:\mathcal{P}\times I\to [0,\infty]$ associated with the distribution $G$ in $\mathcal{P}$ and the measure $\mu$ on $(0,1)$ is defined by\footnote{Recalling that $\mu$ is a Borel measure, the integral (\ref{eq:gds}) is well-defined because the integrand is measurable and $\RG$ is a Borel set.}
\begin{align}
\GDS(F,y;G,\mu)
&= \int_{\RG} S_{\alpha,G^{-1}(\alpha)}(F^{-1}(\alpha),y)\,\dd\mu(\alpha) \notag \\
&= \int_{\RG} (\one\{y<F^{-1}(\alpha)\} - \alpha)(\one\{G^{-1}(\alpha) < F^{-1}(\alpha)\} - \one\{G^{-1}(\alpha) < y\})\,\dd\mu(\alpha), \label{eq:gds}
\end{align}
whenever $F\in\mathcal{P}$ and $y\in I$.

\begin{thm}\label{thm:GDS is equitable}
Suppose that $G\in\mathcal{P}$ and that $\mathcal{P}_0$ is a class of distributions contained in $\mathcal{P}$ such that the expected score $\EEE_G \GDS(F,Y;G,\mu)$ is finite for all $F$ in $\mathcal{P}_0$. Then the generalised diagonal score $\GDS(\twocdots ; G,\mu):\mathcal{P}_0\times I \to[0,\infty]$ is equitable with respect to $G$, with equitability constant $c$ given by
\begin{equation}\label{eq:GDS equitability constant}
c=\int_{\RG}\alpha(1-\alpha)\,\dd\mu(\alpha).
\end{equation}
\end{thm}

When $\dd\mu(\alpha)=\dd\alpha$ and $G$ is continuous, $\GDS(\twocdots;G,\mu)$ recovers the \textit{diagonal score} \citep{ben2018diagonal}
\begin{equation}\label{eq:DS}
\mathrm{DS}(F,y;G)=\int_0^1 S_{\alpha,G^{-1}(\alpha)}(F^{-1}(\alpha),y) \,\dd\alpha,
\end{equation}
while $\mu=\delta_\alpha$ for $\alpha$ in $\RG$ yields the diagonal elementary score $S_{\alpha,G^{-1}(\alpha)}$. They have equitability constants of $1/6$ and $\alpha(1-\alpha)$ respectively. 

A key contribution of this work is the relaxation of the assumption that $G$ is continuous. Equitability is preserved by restricting the domain of integration from $(0,1)$ to $\RG$. For example, suppose that $I=[0,\infty)$ and $\PPP_G(Y=0)=p>0$, so that $G$ has a discontinuity at $0$, as is typical for daily precipitation. Then, by Corollary~\ref{cor:diagonal score equitable}, $S_{\alpha,G^{-1}(\alpha)}$ is not equitable with respect to $G$ if $0<\alpha<p$. Consequently, the diagonal score \eqref{eq:DS}, which integrates over $(0,1)$, is not equitable. If $0$ is the only point at which $G$ is discontinuous, then restricting the domain of integration to $[p,1)$ restores equitability, with a corresponding adjustment to the equitability constant.

\section{Equitable scores for single-valued forecasts}\label{s:real valued equitable}

We now consider equitable scoring functions for single-valued forecasts $x$ in $I$. To the authors' knowledge, only three such functions designed specifically for single-valued forecasts have been proposed in the literature; the case of equitable scores for discretized single-valued forecasts is considered later in Section~\ref{ss:discretize}. Suppose that $G$ is the climatological distribution, and both the forecast $x$ and observation $y$ belong to $I$. The first scoring function, which we denote by $\Grin$, was proposed by \cite{gringorten1965measure} and is defined by
\begin{equation}\label{eq:Grin}
\Grin(x,y)=
\begin{cases}
-\ln(G(x)(1-G(y))) - 1, & x\geq y, \\
-\ln(G(y)(1-G(x))) - 1, & x < y.
\end{cases}
\end{equation}
The second is the so-called `revised LEPS score' \citep{potts1996revised}, which we denote by $\LEPSr$ and is defined by
\begin{equation}\label{eq:LEPSr}
\LEPSr(x,y)=3(1 - |G(x)-G(y)| + G(x)^2-G(x) + G(y)^2 - G(y)) - 1.
\end{equation}
Both are positively oriented scoring functions. The third scoring function, here denoted by $S_1$ and given by
\begin{equation}\label{eq:S1}
S_1(x,y)= (\one\{x>y\}-G(y))^2 - (\one\{x>y\}-G(x))^2,
\end{equation}
is implicit in the work of \citet[Eq.~3]{ben2021verification} via a transformation of the first argument from probability space $[0,1]$ to $I$. It is negatively oriented and was derived from the diagonal score using the same method we use below.

The next proposition shows that, up to equivalence, all three scoring functions are special cases of the generalized diagonal score. To establish this, we first express the generalized diagonal score in the setting of single-valued forecasts $x$.

Given a real number $x$ in $I$, let $F_x$ denote the CDF given by $F_x(z)=\one\{z\geq x\}$, so that $F_x$ represents the distribution concentrated at $x$. In meteorology, such distributions are sometimes called `deterministic'. If we consider the generalised diagonal score restricted to the class $\mathcal{P}_\mathrm{d}$ of deterministic probability distributions, we can identify $F_x$ with the real number $x$ and write
\begin{align}
\GDSd(x,y;G,\mu)
&= \GDS(F_x,y;G,\mu) \notag\\
&= \int_{\RG} S_{\alpha,G^{-1}(\alpha)}(x,y)\,\dd\mu(\alpha) \notag \\
&=  \int_{\RG}  (\one\{y<x\} - \alpha)(\one\{G^{-1}(\alpha) < x\} - \one\{G^{-1}(\alpha) < y\}) \,\dd\mu(\alpha) \label{eq:deterministic GDS}
\end{align}
whenever $x,y\in I$. This representation expresses the score as a weighted average of threshold-based discrepancies indexed by $\alpha$.
Subject to the relevant expectations being finite, this yields equitable scoring functions for single-valued forecasts through different choices of the measure $\mu$. In each case, a lower score is better and a perfect score is $0$.

\begin{table}[t!]
    \renewcommand{\arraystretch}{1.5} 
    \centering \small
    \begin{tabular}{|cl|}
    \hline
        (G1) & $G$ continuous on $\RRR$, $G$ strictly increasing on $\mathcal{D}(G)$ \\ 
        (G2) & $I=[0,\infty)$, $G(0^-)=0<G(0)$, $G$ continuous on $(0,\infty)$, $G$ strictly increasing on $\mathcal{D}(G)$ \\
        (G3) & $G$ satisfies (G1), $0<G(z)<1$ for all $z$ in $\RRR$ \\
        (G4) & $G$ satisfies (G2), $G(z)<1$ for all $z$ in $\RRR$ \\
    \hline
    \end{tabular}
    \caption{Four sets of possible regularity assumptions on the climatological distribution $G$.}
    \label{tab:regularity}
\end{table}

To obtain explicit formulae, we impose regularity conditions on the climatological distribution $G$, summarised in Table~\ref{tab:regularity}. These conditions ensure that the integral in \eqref{eq:deterministic GDS} can be evaluated explicitly, for example by converting indicator functions into intervals of integration, and that certain logarithmic expressions are well-defined. These conditions are satisfied by a wide range of distributions encountered in meteorological applications, ranging from normal distributions (G1, G3) to mixture distributions that arise in variables such as precipitation, streamflow and wind magnitude (G2, G4).

The next proposition derives specific instances of $\GDSd$ under these assumptions and shows how they recover the scoring functions $\Grin$, $\LEPSr$ and $S_1$ up to equivalence. The final part of the proposition provides a flexible class of equitable scoring functions parameterized by a weight function $u$.

\begin{prop}\label{prop:real valued equitable scores}
Suppose that $G\in\mathcal{P}$.
\begin{enumerate}
\item[(i)] Suppose that $G$ satisfies either assumption (G1) or (G2), and the measure $\mu$ has no atoms on $(0,1)$.\footnote{Intuitively, $\mu$ has no atoms if the inclusion or exclusion of endpoints in the interval $[a,b]$ of integration, where $[a,b]\subset(0,1)$, does not change the value of the integral.} Then
\begin{equation}\label{eq:GDSd reduced}
\GDSd(x,y;G,\mu) = \int_{G(y)}^{G(x)}(\one\{x>y\} - \alpha)\,\dd\mu(\alpha) \qquad\forall x,y\in I.
\end{equation}
In particular, the score depends only on the interval between $G(y)$ and $G(x)$ in probability space.
\item[(ii)] If $G$ satisfies either (G1) or (G2) then the scoring function $S_1:I \times I\to[0,1]$ given by (\ref{eq:S1}) is equal to $\GDSd(\twocdots;G,\mu)$ when $\dd\mu(\alpha)=2\,\dd\alpha$, is equitable with respect to $G$, and is equivalent to $\LEPSr$. Its equitability constant is $1/3$ if $G$ satisfies (G1) and $1/3-G(0)^2+2G(0)^3/3$ if  $G$ satisfies (G2).
\item[(iii)] If $G$ satisfies either (G3) or (G4) then the scoring function $S_2: I\times I\to[0,\infty)$, defined by
\begin{equation}\label{eq:S2}
S_2(x,y)=\ln\left|\frac{\one\{x < y\} - G(x)}{\one\{x < y\} - G(y)}\right|, \qquad \forall x,y\in I,
\end{equation}
is equal to $\GDSd(\twocdots;G,\mu)$ when $\dd\mu(\alpha) = (\alpha(1-\alpha))^{-1}\,\dd\alpha$, is equitable with respect to $G$, and is equivalent to $\Grin$. The equitability constant for $S_2$ is $1$ under assumption (G3) and $1-G(0)$ under (G4).
\item[(iv)] Suppose that $u:(0,1)\to\RRR$ is bounded, nonnegative and measurable.\footnote{Measurability is a necessary condition for using $u$ in the integrand. Most functions that practitioners are familiar with are measurable.} Define functions $U$ and $V$ by
\[U(z)=\int_0^z u(\alpha)\,\dd\alpha \quad\text{and}\quad V(z)=\int_0^z \alpha u(\alpha)\,\dd\alpha\]
whenever $0\leq z \leq 1$. If $G$ satisfies either (G1) or (G2) then the scoring function $S_3:I \times I\to[0,\infty)$, defined by
\begin{equation}
S_3(x,y;u)= \one\{x>y\}(U(G(x)) - U(G(y))) - V(G(x)) + V(G(y)), \qquad \forall x,y\in I,
\end{equation}
is equal to $\GDSd(\twocdots;G,\mu)$ when $\dd\mu(\alpha)=u(\alpha)\,\dd\alpha$ and is equitable with respect to $G$. The equitability constant for $S_3$ is $\int_b^1 \alpha(1-\alpha)u(\alpha)\,\dd\alpha$, where $b=0$ under assumption (G1) and $b=G(0)$ under (G2).
\end{enumerate}
\end{prop} 

Note that the scoring functions $S_1$, $S_2$ and $S_3$ are not equitable with respect to $G$ when $G$ has multiple discontinuities. The next proposition therefore provides a score that is equitable with respect to an empirical CDF.

\begin{prop}\label{prop:GDSd ecdf}
Suppose that $G\in\mathcal{P}$ is in the form of a step function, such that
\begin{equation}\label{eq:step func}
G(z)=\sum_{i=1}^N p_i \one\{z_i \leq z\}, \qquad\forall z\in\RRR,
\end{equation}
where $p_i>0$, $\sum_{i=1}^N p_i=1$ and $z_1<z_2<\ldots<z_N$. Let $(w_i)_{i=1}^{N-1}$ be a tuple of nonnegative numbers, and write $P_i=\sum_{k=1}^i p_k$. Then the scoring function $S_4: \RRR\times \RRR\to[0,\infty)$, defined by
\begin{equation}
S_4(x,y) = \sum_{i=1}^{N-1} w_i S_{P_i, z_i}(x,y), \qquad \forall x,y\in\RRR,
\end{equation}
is equal to $\GDSd(\twocdots;G,\mu)$ when $\mu=\sum_{i=1}^{N-1}w_i\delta_{P_i}$, is equitable with respect to $G$, and has equitability constant
\begin{equation}\label{eq:GDSc equitability constant}
\sum_{i=1}^{N-1} w_i P_i(1-P_i).
\end{equation}
\end{prop}

\begin{figure}[t!]
\centering
\includegraphics[width=0.95\textwidth]{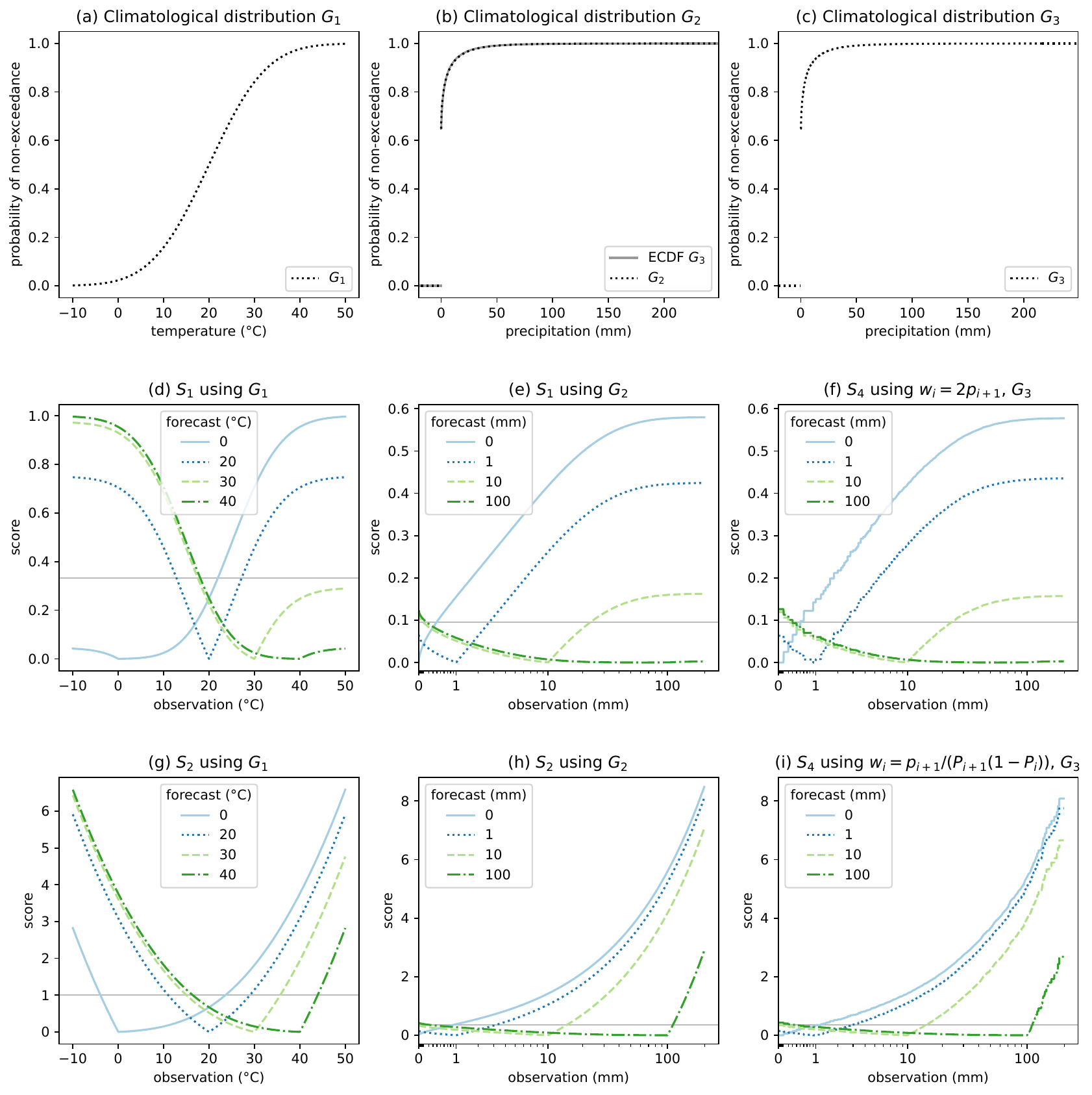}
\caption{(a) The CDF of an idealized climatological distribution $G_1=\cN(20,10^2)$ for temperature. (b) The empirical CDF $G_3$ of the climatological distribution of daily rainfall at Sydney Airport, along with the fitted CDF $G_2$. (c) The empirical CDF $G_3$. Panels (d) to (i): scores of different forecasts against a range of observations when scored using $S_1$, $S_2$ and $S_4$, where each is constructed using one of  $G_1$, $G_2$ and $G_3$ as indicated in each panel title. In the case of $S_4$, the choice of weights $w_i$ is also indicated. The horizontal grey line in each of panels (d) to (i) indicates the equitability constant.}
\label{fig:real valued equitable scores}
\end{figure}

Figure~\ref{fig:real valued equitable scores} shows the scores assigned by $S_1$, $S_2$ and $S_4$ to a range of forecast--observation pairs, such that each scoring function is equitable with respect to at least one of three climatological distributions $G_1$, $G_2$ or $G_3$. The distribution $G_1$ (Fig~\ref{fig:real valued equitable scores}a) is the normal distribution $\cN(20,10^2)$, and models the daily maximum temperature distribution for a location with a wide range of temperatures. The distribution $G_3$ (Fig~\ref{fig:real valued equitable scores}c) is the empirical CDF derived using 34013 quality controlled daily precipitation observations at Sydney Airport, Australia, spanning the period 1 September 1929 to 31 August 2025. The distribution $G_2$ is a parametric fit to this empirical data: we set $G_2(0)=0.648$, which is the proportion of observations that were 0~mm, and then fit the nonzero data to a generalized gamma distribution via the method of moments. As shown in Figure~\ref{fig:real valued equitable scores}b, the fit is very good.  The distribution $G_1$ satisfies assumptions (G1) and (G3), while $G_2$ satisfies (G2) and (G4). The distribution $G_3$ is of the form \eqref{eq:step func}.

In addition to showing the scores assigned to different forecast--observation pairs, Figure~\ref{fig:real valued equitable scores} (d) to (i) also shows the equitability constant for each scoring function as a horizontal grey line. The weights $w_i$ for $S_4$ were chosen to approximate analogous weights induced by the measure $\mu$ of either $S_1$ or $S_2$. For each of these scores, whenever a forecast lies in the right-tail of the climatological distribution, underprediction is penalised more harshly than a similarly sized (in absolute terms) overprediction. Thus, when the observation is likely in the tail of the climatological distribution, these scores encourage forecasts further into the tail. We make this observation precise in Section~\ref{s:optimal}.

\section{Equitable scores for multicategorical forecasts}\label{s:categorical equitable}

We now turn to equitable scoring functions for multicategorical forecasts, which comprise the majority of equitable scores studied in the literature. Let $N$ denote the number of categories, and let $[N] := \{1,2,\ldots,N\}$, where each $i$ in $[N]$ represents the $i$th category. We assume the categories are ordinal, so that their ordering carries meaning.

In some applications, the categories correspond to disjoint intervals $C_1,\ldots,C_N$ that partition $I$ (for example, when a continuous variable is discretized). In such cases, category $C_i$ may be identified with $i \in [N]$ when convenient.

Suppose that $Y$ is a random variable taking values in $[N]$ with climatological distribution $G$ whose mass is concentrated on $[N]$. For each $k$ in $[N]$, define $p_k$ to be the climatological probability of category $k$ and $P_k$ to be the climatological cumulative probability for category $k$. Unless stated otherwise, we assume that $p_k > 0$ for all $k$. Thus
\begin{align*}
p_k &= \PPP_G(Y = k), \qquad & \forall k \in [N], \\
P_k &= \PPP_G(Y \le k) = \sum_{i=1}^k p_i, \qquad &  \forall k \in [N], \\
\end{align*}
and $\sum_{k=1}^N p_k = 1$. The CDF $G$ can then be written as
\begin{equation}\label{eq:G categorical}
G(z)=\sum_{k=1}^N p_k \one\{k \leq z\}, \qquad \forall z \in\RRR,
\end{equation}
so that $G$ is a step function of the form \eqref{eq:step func} with jump points $z_k$ occurring when $z_k = k$. Thus $G$ is completely determined by the climatological probabilities $(p_1,\ldots,p_N)$.

It follows from Proposition~\ref{prop:GDSd ecdf} that, for such a step function $G$, the generalized diagonal score $\GDSd$ with a discrete measure $\mu$ supported on $\{P_1,\ldots,P_{N-1}\}$ yields an equitable scoring function. Restricting the forecast and observation domains to $[N]$, we therefore obtain an equitable scoring function on $[N]\times[N]$.

We summarise this construction as follows. Let $\boldP$ denote the vector of climatological cumulative probabilities $(P_1,\ldots,P_N)$ and let $\boldw$ denote a vector $(w_1,\ldots,w_{N-1})$ of nonnegative weights. Define the scoring function $\GDSc(\,\cdot\,,\,\cdot\,;\boldP,\boldw):[N]\times[N] \to [0,\infty)$ by
\begin{align}
\GDSc(i,j; \boldP, \boldw)
&= \GDSd\!\left(i,j;G,\sum_{k=1}^{N-1} w_k \delta_{P_k}\right) \notag \\
&= \sum_{k=1}^{N-1} w_k S_{P_k,k}(i,j),
\label{eq:GDSc}
\end{align}
for all $i$ and $j$ in $[N]$. Thus $\GDSc$ aggregates binary threshold-based penalties across category boundaries, with weights $w_k$ controlling the relative importance of errors at different thresholds. The preceding construction yields the following corollary.

\begin{cor}\label{cor:GDSc}
Suppose that $G$ is a probability distribution on the space $[N]$ of categories with cumulative probability vector $\boldP=(P_1,\ldots,P_N)$. If $\boldw=(w_1,\ldots,w_{N-1})$ is a vector of nonnegative weights, then $\GDSc(\twocdots; \boldP, \boldw):[N]\times[N]\to[0,\infty)$ is equitable with respect to $G$ with equitability constant given by \eqref{eq:GDSc equitability constant}.
\end{cor}

The scoring function $\GDSc$ can be represented by a scoring matrix $M$, whose $(i,j)$th entry $m_{i,j}$ gives the penalty incurred when category $i$ is forecast and category $j$ is observed. Explicitly,
\begin{equation}\label{eq:S_c matrix entries}
m_{i,j} =\GDSc(i,j;\boldP,\boldw)=
\begin{cases}
0, & 1\leq i=j \leq N, \\
\sum_{k=i}^{j-1} w_k P_k , & 1\leq i<j \leq N, \\
\sum_{k=j}^{i-1} w_k (1-P_k) , & 1\leq j<i \leq N,
\end{cases}
\end{equation}
Scores here are interpreted as a penalty, so that a lower score is better and a score of zero indicates the forecast and observed categories were identical. Note that $\GDSc(\twocdots;\boldP,\boldw)$, as with the Gandin--Murphy scores, penalises larger category discrepancies more heavily than smaller ones.

We examine the 2- and 3-category variants of $\GDSc$. In the 2-category case, the scoring matrix $M$ is given by
\[M
=\begin{pmatrix}
0 & w_1 P_1 \\
w_1(1-P_1) & 0 
\end{pmatrix}.
\]
An equivalent positively oriented scoring matrix $M'$, obtained via $m'_{i,j}=a m_{i,j}+h(j)$ where $a=-(w_1 P_1(1-P_1))^{-1}$ and
\[
h(j)=
\begin{cases}
(1-P_1)/P_1, & j=1, \\
P_1/(1-P_1), & j=2,
\end{cases}
\]
is given by
\[
M'
=\left(\begin{matrix}
\ds\frac{1-P_1}{P_1} & -1 \\
-1 & \ds\frac{P_1}{1-P_1}
\end{matrix}\right)
=\left(\begin{matrix}
\ds\frac{p_2}{p_1} & -1 \\
-1 & \ds\frac{p_1}{p_2}
\end{matrix}\right).
\]
The scoring matrix $M'$ is the \citet{peirce1884numerical} skill score \citep[Table~2]{rodwell2010new} and, along with equivalent scores, is noted by \citet{gandin1992equitable} as being `the only performance measures in the two-event situation that satisfy the basic conditions of equitability'.


In the 3-category case, the scoring matrix $M$ obtained from $\GDSc$ equals
\begin{equation}\label{eq:GDS scoring matrix 3-cat}
\begin{pmatrix}
0 & w_1 P_1 &  w_1 P_1+w_2 P_2 \\
w_1(1- P_1) & 0 & w_2 P_2  \\
w_1(1-P_1) + w_2(1-P_2) & w_2(1-P_2) & 0
\end{pmatrix}.
\end{equation}
The 3-category Gerrity score with climatological probabilities $(p_1,p_2,p_3)$ has scoring matrix
\begin{equation}\label{eq:Gerrity 3-cat}
\frac{1}{2}\begin{pmatrix}
\ds\frac{1-p_1}{p_1} + \frac{p_3}{1-p_3} & \ds\frac{p_3}{1-p_3}-1 & -2 \\
\ds\frac{p_3}{1-p_3} -1 & \ds\frac{p_1}{1-p_1} + \frac{p_3}{1-p_3} &  \ds\frac{p_1}{1-p_1} -1 \\
-2 & \ds\frac{p_1}{1-p_1} - 1 & \ds\frac{p_1}{1-p_1} + \ds\frac{1-p_3}{p_3}
\end{pmatrix},
\end{equation}
while the SEEPS score with climatological probabilities $(p_1,p_2,p_3)$ has scoring matrix
\begin{equation}\label{eq:SEEPS}
\frac{1}{2}\begin{pmatrix}
0 & \ds\frac{1}{1-p_1} & \ds\frac{1}{p_3}+\frac{1}{1-p_1} \\
\ds\frac{1}{p_1} & 0 &  \ds\frac{1}{p_3} \\
 \ds\frac{1}{p_1}+\frac{1}{1-p_3} & \ds\frac{1}{1-p_3} & 0
\end{pmatrix}.
\end{equation}
A consequence of the next proposition is that, subject to appropriate choices of weights $w_i$, these three scores are equivalent.

\begin{prop}\label{prop:gerrity SEEPS GDS equivalence}
Suppose that $N$ ordinal categories $[N]$ have corresponding positive climatological probabilities $(p_1,\ldots\,p_N)$ and cumulative probability vector $\boldP=(P_1,\ldots,P_N)$. Define $\boldw=(w_1,\dots,w_{N-1})$ by $w_k =  ((N-1) P_k(1-P_k))^{-1}$ whenever $1\leq k \leq N-1$.
\begin{itemize}
\item[(i)] The SEEPS score with scoring matrix (\ref{eq:SEEPS}) equals $\GDSc(\,\cdot\,,\,\cdot\,; \boldP, \boldw)$ when $N=3$.
\item[(ii)] The positively oriented Gerrity score for 3 ordinal categories with corresponding climatological probabilities $(p_1,p_2,p_3)$ is equivalent to the negatively oriented SEEPS score with the same categorical climatological probabilities.
\item[(iii)] The positively-oriented Gerrity score for $N$ ordinal categories $[N]$ is equivalent to the negatively oriented scoring function $\GDSc(\,\cdot\,,\,\cdot\,; \boldP, \boldw)$.
\item[(iv)] The equitability constant of $\GDSc(\,\cdot\,,\,\cdot\,; \boldP, \boldw)$ is $1$.
\end{itemize}
\end{prop}

Part (ii) of the proposition was essentially recognised by \citet[Section~5.2.1]{rodwell2010new}. Our proof shows that a forecaster's SEEPS score is the difference between the Gerrity score awarded for a perfect prediction and the Gerrity score awarded to the forecaster's prediction. That is, the SEEPS score can be interpreted as the forecaster's regret, relative to an omniscient oracle, when rewards are allocated by the Gerrity score (c.f. Section~2.3 of \citealt{ehm2016quantiles}). The same relationship exists between the Gerrity score and the instance of $\GDSc$ in part (iii).

Other equitable scores for categorical forecasts in the literature are equivalent to specific instances of $\GDSc$. For example, the equitable adjustment of the Heidke score developed by \citet[Table~4]{barnston1992correspondence} is equivalent to $\GDSc(\,\cdot\,,\,\cdot\,; \boldP, \boldw)$ where $w_k=1$ and the cumulative probabilities satisfy $P_k=k/N$. The proof of equivalence follows the same method of proof as Proposition~\ref{prop:gerrity SEEPS GDS equivalence}(ii).

All the scoring functions discussed in this section belong to the class of categorical equitable scoring functions developed by \citet{gandin1992equitable}. However, not every Gandin--Murphy equitable score is equivalent to the generalised diagonal score for suitable measure $\mu$. To give an example, we consider the categorical version of $\LEPSr$ proposed by \citet{potts1996revised}, called LEPSCAT by \citet{Livezey2003Categorical}, and denoted here by $\LEPSc$.

The $\LEPSc$ score is derived from $\LEPSr$ by partitioning the outcome space $\RRR$ for a continuous random variable $Y$ into $N$ ordinal categories $(C_1,C_2,\ldots,C_N)$. For continuous $G$,  $\LEPSc:\{C_1,\ldots,C_N\}\times\{C_1,\ldots,C_N\}\to\RRR$ is defined by
\begin{equation}\label{eq:LEPSCAT}
\LEPSc(C_i,C_j) = \EEE[\LEPSr(X,Y) \mid X\sim G, Y\sim G, X\in C_i, Y\in C_j],
\end{equation}
where the random variables $X$ and $Y$ are independent.\footnote{\citet{potts1996revised} write: ``The scores for forecasts issued in discrete categories are obtained by calculating the expected score for each combination of observed and forecast categories.'' We interpret this as the expectation with respect to the climatological distribution $G$, and further impose the assumption that $G$ is continuous so that all two scoring matrices they published are recovered.} The $\LEPSc$ score has positive orientation. \citet{potts1996revised} did not provide a general formula for $\LEPSc$, nor a proof that $\LEPSc$ is equitable, and only supplied scoring matrices for the 3- and 5-category cases of equiprobable categories. This lack of an explicit formula for the scoring matrix provides some explanation for why \citet{Livezey2003Categorical} recommends the Gerrity score over $\LEPSc$. We rectify this situation with the following proposition. In the following, let $P_i$ denote the cumulative climatological probabilities of categories $C_1$ through to $C_i$, and define $P_0=0$.

\begin{prop}\label{prop:LEPSCAT}
Suppose that the climatological distribution $G$ in $\mathcal{P}$ is continuous and that $q_{i,j}'$ denotes the $(i,j)$th entry of the $\LEPSc$ scoring matrix. Then
\begin{align}
q_{i,i}' &= 2(P_i^2 +P_{i-1} P_i+P_{i-1}^2-2 P_i-P_{i-1}+1),\quad & 1\leq i \leq N, \label{eq:LEPSCAT entry 1}\\
q_{i,j}' &= 2 - 3(P_{i-1}+P_i) + 3 P_{i-1} P_i + 3 P_{j-1} P_j + (P_i-P_{i-1})^2 + (P_j-P_{j-1})^2,\quad & 1\leq j < i \leq N,  \label{eq:LEPSCAT entry 2} \\
q_{i,j}' &= q_{j,i}', & 1\leq i < j \leq N.  \label{eq:LEPSCAT entry 3}
\end{align}
Moreover, $\LEPSc$ is equitable with respect to $G$ with equitability constant $0$.
\end{prop}

The 3-category $\LEPSc$ scoring matrix with equiprobable climatological probabilities $(p_i=1/3)$ is given by
\begin{equation}\label{eq:LEPSCAT 3 equi}
\frac{1}{9}
\begin{pmatrix}
8 & -1 & -7 \\
-1 & 2 & -1 \\
-7 & -1 & 8
\end{pmatrix}
\end{equation}
\citep[Table~1]{potts1996revised}. This is not equivalent to a score in $\GDSc$ form.

\begin{prop}\label{prop:gandin murphy not GDS}
The scoring matrix (\ref{eq:LEPSCAT 3 equi}) and the scoring matrices of Equations~(27) and (32) of \citet{gandin1992equitable} are not equivalent to any scoring matrix $M$ of the form (\ref{eq:S_c matrix entries}). Consequently, up to equivalence, the family of categorical generalized diagonal scores $\GDSc$ is a proper subset of the family of Gandin--Murphy equitable scores.
\end{prop}

Notwithstanding Proposition~\ref{prop:gandin murphy not GDS}, the generalized diagonal score encompasses many equitable scores of practical interest.

\subsection{Equitable scores for discretized real-valued forecasts and observations}\label{ss:discretize}

We now consider the case where forecasts and observations are real-valued but are subsequently discretized (or categorized) prior to being scored. Prominent examples include the verification of precipitation forecasts from weather prediction models using SEEPS \citep{rodwell2010new}, and the discretisation of real-valued temperatures into terciles or quintiles in seasonal forecasting \citep{barnston1992correspondence, potts1996revised}.

We begin by partitioning the outcome space $I$ into $N$ categories $C_1,\ldots,C_N$ using thresholds $\theta_1,\ldots,\theta_{N-1}$ in $I$ satisfying $\theta_1<\cdots<\theta_{N-1}$ as follows:
\begin{align*}
C_1&=\{z\in I:z\leq \theta_1\}, \quad & \\
C_i&=(\theta_{i-1},\theta_i], \quad & 1<i<N, \\
C_N&=\{z\in I:z>\theta_{N-1}\}. \quad &
\end{align*}

\begin{prop}\label{prop:equitable discretize}
Suppose that $G\in\mathcal{P}$, and that $\theta_i\in \mathcal{D}(G)$, $G^{-1}(G(\theta_i))=\theta_i$ and $w_i>0$ whenever $1\leq i \leq N-1$. Define the scoring function $S_5:I \times I\to[0,\infty)$ by
\begin{equation}\label{eq:equitable discretize}
S_5(x,y)=\sum_{i=1}^{N-1}w_i S_{G(\theta_i), \theta_i}(x,y)\,\qquad \forall x,y\in I.
\end{equation}
Then $S_5$ is equal to $\GDSd(\twocdots; G, \mu)$ when $\mu=\sum_{i=1}^{N-1}w_i\delta_{G(\theta_i)}$, is equitable with respect to $G$ and has equitability constant
\begin{equation}
\sum_{i=1}^{N-1} w_i G(\theta_i)(1-G(\theta_i)).
\end{equation}
Moreover, if $x,x'\in C_j$ and $y,y'\in C_k$ then $S_5(x,y)=S_5(x',y')$.
\end{prop}

The final part of the proposition shows that the score $S_5(x,y)$ depends only on the categories to which the forecast $x$ and observation $y$ belong, and not on their exact values. Consequently, $S_5$ can be interpreted as an equitable score in the purely categorical setting. In this case, its scoring matrix has entries of the form \eqref{eq:S_c matrix entries} with $P_k = G(\theta_k)$. 

The condition that $\theta_i \in \mathcal{D}(G)$ and $G^{-1}(G(\theta_i))=\theta_i$ is satisfied if and only if $\theta_i = \min\{ z \in \mathcal{D}(G) : G(z) = G(\theta_i) \}$. Consequently, if $\theta_i$ lies in an interval on which $G$ is constant, then it must coincide with the left endpoint of that interval to satisfy the hypothesis of the proposition.

We now illustrate this construction using the SEEPS score, adopting the parameter settings recommended by \citet{rodwell2010new} for the verification of 24-hour precipitation forecasts. Here, both forecasts and observations take values in the interval $I=[0,\infty)$, in units of millimetres. Let $G$ denote the climatological CDF at a given location and suppose that $G$ satisfies condition (G2) of Table~\ref{tab:regularity}.
We partition $I$ into three categories $C_1$ (`dry'), $C_2$ (`light') and $C_3$ (`heavy') defined by
\begin{equation}\label{eq:SEEPS categories}
C_1 = [0,\theta_1], \,\, C_2 = (\theta_1,\theta_2],  \,\, C_3 = (\theta_2,\infty) \quad\text{where}\quad \theta_1=0.2, \,\, p_1 = G(\theta_1),\,\, \theta_2 = G^{-1}((2+p_1)/3).
\end{equation}
Thus $p_1$ is the climatological probability of the dry category, and the probability of $C_2$ is twice that of $C_3$. Then with weights $w_i$ chosen such that $w_i=(2 G(\theta_i) (1-G(\theta_i)))^{-1}$, the score $S_5$ coincides with the SEEPS score applied to real-valued forecasts and observations. The same construction works for an empirical CDF $G$, provided that 0.2 is a jump point for $G$ and that $2(1-p_1)/3 \in\RG$.

With this example, we have now shown that every scoring function in Table~\ref{tab:scores in GDS form} is equivalent to a score in GDS form.

\section{Optimal forecasts of equitable scores}\label{s:optimal}

In this section we discuss optimal forecast strategies when being assessed with equitable scores. In particular, suppose that the forecaster's beliefs about the outcome $Y$ are expressed as a predictive distribution $F$ in $\mathcal{P}$. We ask the following question: Which forecasts $f$ in the domain $\mathcal{F}$ of possible forecasts minimize the expected score $\EEE_F S(f,Y)$, where $S:\mathcal{F}\times\mathcal{O}\to[0,\infty]$ is a scoring function in $\GDS$ form. For brevity, a forecast that minimizes the expected score is said to be optimal. Such forecasts are optimal for any equivalent scoring function.

\subsection{Optimal predictive distributions for $\GDS$}\label{ss:optimal GDS}

First we fix $G\in\mathcal{P}$ and the measure $\mu$. Let $\mathcal{P}_0$ denote a class of distributions such that $\EEE_F \GDS(\tilde{F},Y;G,\mu)$ exists and is finite whenever $F,\tilde{F}\in\mathcal{P}_0$. For example, if $\dd\mu(\alpha)=\dd\alpha$, which recovers the diagonal score, then one may take $\mathcal{P}_0$ to be the class of distributions in $\mathcal{P}$ with finite first moment.

\begin{defn}
A scoring rule $S$ is said to be \textit{proper} relative to a class of distributions $\mathcal{P}_0$ if 
\begin{equation}\label{eq:proper}
\EEE_F\GDS(F,Y;G,\mu) \leq \EEE_F\GDS(\tilde{F},Y;G,\mu)
\end{equation}
whenever $F,\tilde{F}\in\mathcal{P}_0$. It is said to be \textit{strictly proper} if (\ref{eq:proper}) holds with equality if and only if $\tilde{F}=F$.
\end{defn}

Thus strictly proper scoring rules encourage honest forecasting \citep{gneiting2007strictly} because the forecaster's expected score is uniquely optimized by issuing the forecast that matches their true probabilistic beliefs.

\begin{thm}\label{thm:proper not strict}
Suppose that $\mathcal{P}_0$ denotes a class of distributions such that the expected score $\EEE_F \GDS(\tilde{F},Y;G,\mu)$ is finite whenever $F,\tilde{F}\in\mathcal{P}_0$.
\begin{enumerate}
\item[(i)] The scoring rule $\GDS(\twocdots;G,\mu)$ is proper relative to $\mathcal{P}_0$.
\item[(ii)] If $\mathcal{P}_0$ contains $G$ and at least one other distinct distribution $F$, then
\begin{equation}\label{eq:proper not strict}
\EEE_G\GDS(G,Y;G,\mu) = \EEE_G\GDS(F,Y;G,\mu)
\end{equation}
and hence $\GDS(\twocdots;G,\mu)$ is not strictly proper relative to $\mathcal{P}_0$.
\item[(iii)] Suppose that $\mathcal{P}_0$ contains two distinct distributions $F$ and $\tilde{F}$ that satisfy the following property:
\begin{equation}\label{eq:same sign condition}
\text{For all $\alpha$ in $(0,1)$,}\quad  F^{-1}(\alpha) \leq G^{-1}(\alpha) \quad \text{if and only if}\quad \tilde{F}^{-1}(\alpha) \leq G^{-1}(\alpha).
\end{equation}
Then
\[
\EEE_F\GDS(F,Y;G,\mu) = \EEE_F\GDS(\tilde{F},Y;G,\mu)
\]
and hence $\GDS(\twocdots;G,\mu)$ is not strictly proper relative to $\mathcal{P}_0$.
\end{enumerate}
\end{thm}

Note that Equation~\eqref{eq:proper not strict} is an immediate consequence of equitability: all predictive distributions are optimizers of the expected generalized diagonal score when the data-generating distribution is the climatological distribution $G$. Hence, whenever $G$ is one possible forecast among others, the generalized diagonal score cannot be strictly proper. In particular, this answers in the negative the question posed by \citet{ben2018diagonal} as to whether the diagonal score is strictly proper. 

\begin{figure}[t!]
\centering
\includegraphics[width=0.9\textwidth]{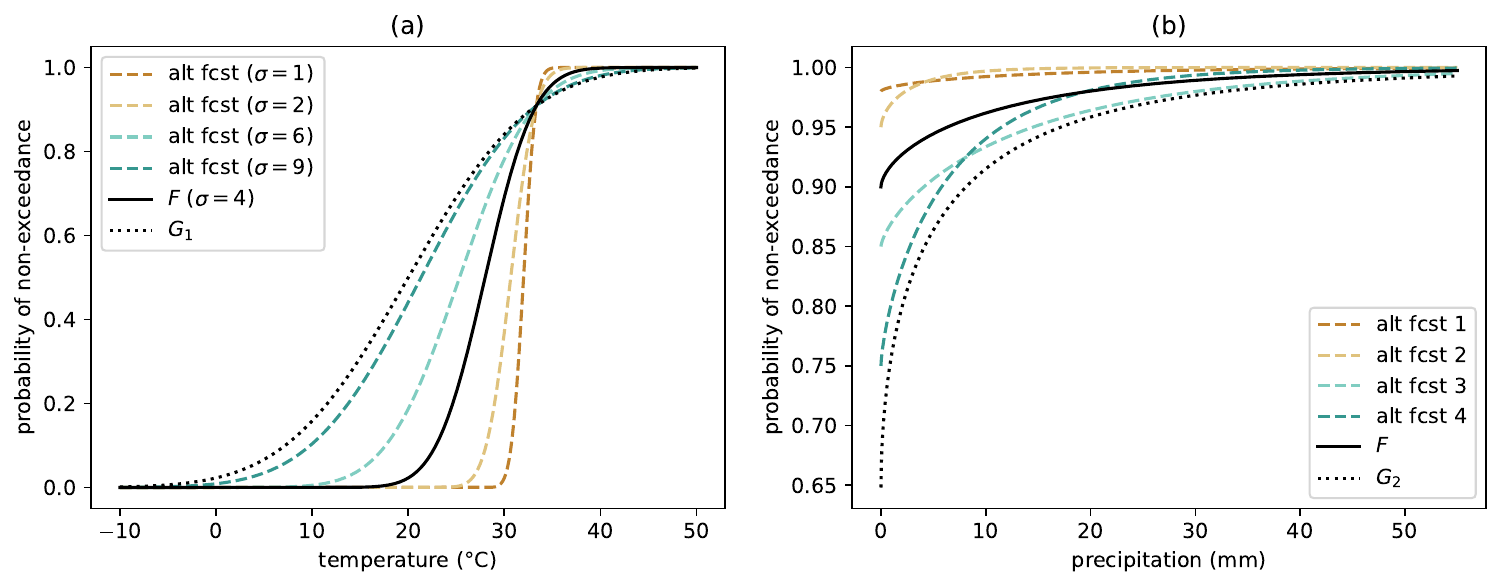}
\caption{In each panel, the CDFs of a  climatological distribution $G_i$, a predictive distribution $F$ (encoding a forecaster's true beliefs) and four alternative forecasts are graphed. The alternative forecasts and $F$ all minimize the forecaster's expected diagonal score, because the alternative forecasts always lie on the same side of the climatological CDF as $F$. In panel (b), all CDFs are for nonnegative random variables and have a discontinuity at 0.}
\label{fig:normal crossing point}
\end{figure}

Part (iii) highlights that the GDS is sensitive to whether each predictive quantile lies above or below the corresponding climatological quantiles, but is insensitive to how far above or below each lies. To illustrate the consequences of this, consider the score $\GDS(\twocdots;G_1,\mu)$, where $G_1=\cN(20, 10^2)$ (c.f. Fig~\ref{fig:real valued equitable scores}a) and where $\dd\mu(\alpha)=\dd\alpha$ (i.e., the measure that yields the diagonal score). Suppose that the distribution $F=\cN(28,4^2)$ best represents the forecaster's uncertainty about the daily maximum temperature for a particular forecast case. Figure~\ref{fig:normal crossing point}a shows graphs of the CDFs $G_1$ and $F$, along with four alternative forecasts $\tilde{F}$, each having a normal distribution with different standard deviation $\sigma$ ($\sigma\in\{1,2,6,9\}$). Since the graph of each alternative forecast $\tilde{F}$ lies on the same side of $G_1$ as $F$, and all the CDFs are continuous and strictly increasing, condition (\ref{eq:same sign condition}) holds for each $\tilde{F}$. It follows from Theorem~\ref{thm:proper not strict}(iii) that $F$ and each of the alternative forecasts attain the minimum expected score and are thus, from the forecaster's perspective, optimal for attaining the best diagonal score. Figure~\ref{fig:normal crossing point}b illustrates the same point, this time using CDFs of nonnegative random variables with discontinuities at $0$. Here $G_2$ is the fitted climatological CDF for precipitation at Sydney airport (c.f. Fig~\ref{fig:real valued equitable scores}b). The CDFs of the forecaster's predictive distribution $F$ and four alternative forecasts $\tilde{F}$ also satisfy condition (\ref{eq:same sign condition}) and thus have the same expected diagonal score. These examples illustrate that forecasts can be substantively misspecified relative to the forecaster's judgment yet still score optimally in expectation.

\citet{jolliffe2008proper} showed that, in the context of probabilistic forecasts where the outcome is binary, there are no scoring rules that are both equitable and proper, under the assumption that the scoring rule is sufficiently smooth. By relaxing this assumption, the next theorem gives a new family of scoring rules for probabilistic forecasts with categorical outcomes that are simultaneously proper and equitable. These are derived from the $\GDS$ for judicious choice of measure $\mu$, and the properties of equitability and propriety follow from Theorems~\ref{thm:GDS is equitable} and \ref{thm:proper not strict}. We present the general case first, followed by the special case of binary outcomes.

Recall that $[N]$ represents the set of $N$ ordinal categories. Let $\mathcal{P}_\mathrm{c}$ denote the space of probability distributions on $[N]$, so that each distribution in $\mathcal{P}_\mathrm{c}$ can be represented by a vector $(r_1,r_2,\ldots,r_N)$ of nonnegative components that sum to $1$, where $r_i$ is the probability that the random variable equals category $i$. Given a  climatological distribution $(p_1,\ldots, p_N)$ in  $\mathcal{P}_\mathrm{c}$, and nonnegative weights $w_k$ ($k=1,\ldots,N-1$) define the scoring rule $S:\mathcal{P}_\mathrm{c}\times[N]\to[0,\infty)$ by
\begin{equation}\label{eq:proper equitable for cats}
S((r_1,\ldots,r_N), i) = \sum_{k=1}^{N-1} w_k(\one\{R_i<P_k\}-P_k)(\one\{R_k<P_k\} - \one\{k<i\})
\end{equation}
whenever $(r_1,\ldots,r_N)\in\mathcal{P}_\mathrm{c}$ and $i\in[N]$. Here $R_k=\sum_{j=1}^k r_j$ and $P_k=\sum_{j=1}^k p_j$ denote the cumulative probabilities. 

\begin{cor}\label{cor:proper equitable for cats}
Suppose that the climatological distribution $(p_1,\ldots, p_N)\in \mathcal{P}_\mathrm{c}$ satisfies $p_i>0$ whenever $i\in[N]$. Then the scoring rule $S:\mathcal{P}_\mathrm{c}\times[N]\to[0,\infty)$ given by (\ref{eq:proper equitable for cats}) is equitable with respect to that climatological distribution and is proper relative to $\mathcal{P}_\mathrm{c}$. The equitability constant is given by \eqref{eq:GDSc equitability constant}.
\end{cor}

In the binary case, the outcome space is often represented by $\{0,1\}$, where outcome $1$ is an event. Identifying category $1$ with $y=0$ and category $2$ with $y=1$, and writing $p$ and $r$ for the climatological and predictive probabilities of an event, respectively, the scoring rule (\ref{eq:proper equitable for cats}) with $w_1=1$ reduces to
\begin{equation}\label{eq:proper equitable for binary}
S(r,y)=p\one\{r>p, y=0\} + (1-p)\one\{r\leq p, y=1\}.
\end{equation}
This proper equitable scoring rule assigns asymmetric penalties for any probability predictions $r$ that are less confident of the observed outcome $y$ than climatology.

\subsection{Optimal single-valued forecasts for $\GDSd$ and categorical forecasts for $\GDSc$}

The condition \eqref{eq:same sign condition} of Theorem~\ref{thm:proper not strict}, illustrated in Figure~\ref{fig:normal crossing point}, is closely related to the concept of `crossing points', i.e., points where the difference between the predictive and climatological CDFs changes sign from negative to positive \citep{ben2021verification}. This concept arises naturally when studying optimal forecasts for equitable scores of single-valued and categorical forecasts that are expressible in GDS form.

\begin{defn}\label{def:crossing point}
Suppose $F,G\in\mathcal{P}$ and $I_0$ is an interval. We say that a point $\chi$ in $I_0$ is a \textit{crossing point} for $F$ on $I_0$ with respect to $G$ if $F(z)\leq G(z)$  whenever $z\in(-\infty,\chi)\cap I_0$ and $F(z)\geq G(z)$  whenever $z\in[\chi,\infty)\cap I_0$.
\end{defn}

\begin{figure}[t!]
\centering
\includegraphics[width=0.98\textwidth]{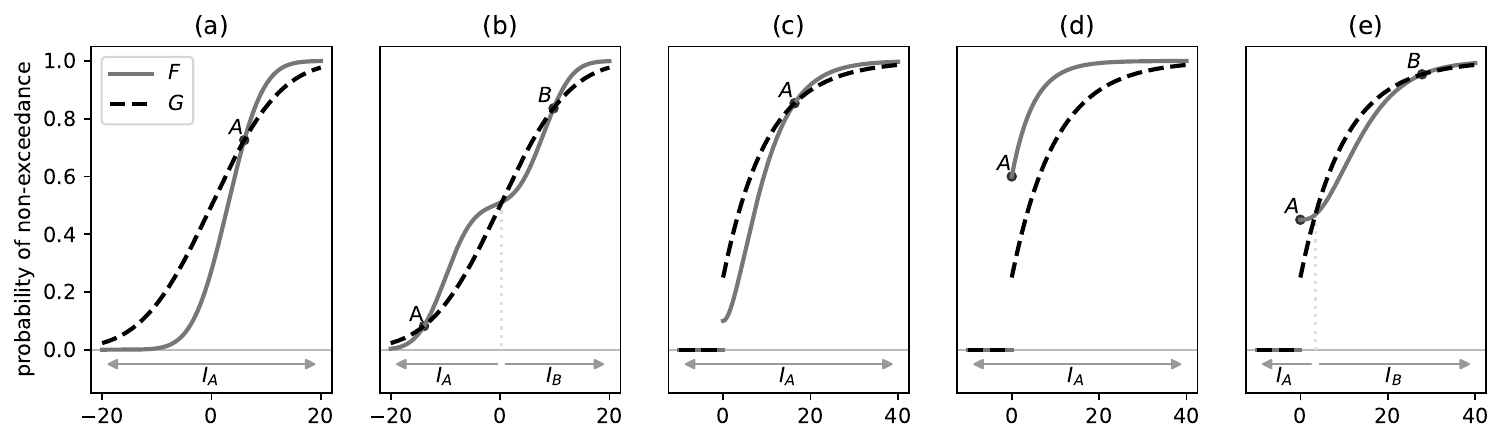}
\caption{Five typical configurations between a climatological CDF $G$ (dashed curve) and predictive CDF $F$ (solid curve), where both are right-continuous and, in panels (c) to (e), have a discontinuity at $0$. In each case, the abscissa $\chi_A$ is a crossing point for $F$ on the interval $I_A$ with respect to $G$, with interval $I_A$ either the entire real line or a ray as depicted using arrows. In cases (b) and (e), the abscissa $\chi_B$ is a crossing point for $F$ on the interval $I_B$ with respect to $G$.}
\label{fig:optimal crossing point}
\end{figure}

Figure~\ref{fig:optimal crossing point} illustrates five typical configurations between a climatological distribution $G$ (dashed curve) and predictive distribution $F$ (solid curve) that occur in meteorological applications (c.f. \citealt[Fig.~6]{ben2021verification}). In each case, the point $A$ has coordinates $(\chi_A,F(\chi_A))$, and the point $\chi_A$ is a crossing point for $F$ on the interval $I_A$ with respect to $G$. Similarly, $B$ has coordinates $(\chi_B,F(\chi_B))$, and  $\chi_B$ is a crossing point for $F$ on the interval $I_B$ with respect to $G$. Note that the definition allows endpoints of closed intervals to be crossing points. For example, in Figure~\ref{fig:optimal crossing point}d, $\chi_A=0$ is a crossing point for $F$ on the interval $[0,20]$.

\begin{thm}\label{thm:optimal}
Suppose that $F,G\in\mathcal{P}$, that $a<b$ and that the expected score $\EEE_F \GDSd(x,Y;G,\mu)$ is finite for every $x$ in $[a,b]$. If $\chi$ is a crossing point for $F$ on $[a,b]$ with respect to $G$, then
\begin{equation}\label{eq:optimal expectation}
\EEE_F \GDSd(\chi,Y;G,\mu) \leq \EEE_F \GDSd(x,Y;G,\mu)
\end{equation}
for all $x$ in $[a,b]$.
\end{thm}

Essentially, if $\chi$ is a crossing point on $[a,b]$ then $\chi$ is an optimal forecast for $\GDSd$ on $[a,b]$. We state the theorem for finite closed intervals $[a,b]$ because that is sufficient for the proof. Nonetheless, it also extends to intervals that are closed relative to $I$, including unbounded intervals, by applying the theorem to suitable finite closed subintervals. For example, we obtain the following corollary.

\begin{cor}\label{cor:crossing point}
Suppose $F,G\in\mathcal{P}$ and that $\EEE_F \GDSd(x,Y;G,\mu)$ is finite for every $x$ in $I$. If $\chi$ is a crossing point for $F$ on $I$ with respect to $G$ then the mapping $x\mapsto\EEE_F \GDSd(x,Y;G,\mu):I\to[0,\infty)$ attains a global minimum at $\chi$.
\end{cor}

In plain language, if there is a unique crossing point then it is the optimal forecast for any score in $\GDSd$ form. In the case where real-valued forecasts and observations are discretized, the corollary becomes the following: if there is a unique crossing point $\chi$, then any single-valued forecast in the same category as $\chi$ is an optimal forecast for the scoring function $S_5$ \eqref{eq:equitable discretize}. In particular, to optimize SEEPS, any forecast lying in the same category as the crossing point is optimal.

Theorem~\ref{thm:optimal} and Corollary~\ref{cor:crossing point} also apply to scores in $\GDSc$ form \eqref{eq:GDSc}, as these are special cases of $\GDSd$. In the categorical context, the predictive and climatological distributions can be identified with their cumulative probability vectors $(R_i)_{i=1}^N$ and $(P_i)_{i=1}^N$. For any integers $i$ and $j$ satisfying $1\leq i<j\leq N$, let $[i:j]$ denote the set $\{i,i+1,\dots,j\}$. In this setting $\chi\in[i:j]$ is a crossing point on $[i:j]$ if $R_k\leq P_k$ for all $k<\chi$ in $[i:j]$ and $R_k\geq P_k$ for all $k\geq \chi$ in $[i:j]$. If $\chi$ is a crossing point on $[i:j]$ then category $\chi$ is an optimal forecast among the categories in $[i:j]$ for any scoring function in $\GDSc$ form, including the SEEPS score, and by equivalence also for the Gerrity score. If $\chi$ is the unique crossing point on $[N]$ then $\chi$ is an optimal forecast category for $\GDSc$.

The preceding results hold provided that the expectation $\EEE_F \GDSd(x,Y;G,\mu)$ is finite. The next proposition addresses this issue for the specific instances of $\GDSd$ considered so far. A sufficient condition for the finiteness of the expected $S_2$ score is a tail dominance property. When this condition holds, it indicates that the predictive distribution $F$ is, in a particular sense, `sharper' or `more concentrated' than the climatological distribution $G$. This property holds in most meteorological applications.

\begin{defn}\label{def:global sharp}
Suppose that $F,G\in\mathcal{P}$. We say that $F$ is \textit{globally sharper} than $G$ on $I$ if there exist constants $a_0$ and $b_0$ such that $F(z)\leq G(z)$ whenever $z\leq a_0$ and $F(z)\geq G(z)$ whenever $z\geq b_0$.
\end{defn}

\begin{prop}\label{prop:optimal fcst}
Suppose that $F, G \in\mathcal{P}$.
\begin{enumerate}
\item[(i)] If $S$ is one of $S_1$, $S_3$, $S_4$, $S_5$ or $\GDSc$ then $\EEE_F S(x,Y)$ is finite for every $x$ in $I$.
\item[(ii)] If $0<G(z)<1$ whenever $z\in I$ and $F$ is globally sharper than $G$, then $\EEE_F S_2(x,Y)$ is finite for every $x$ in $I$.
\end{enumerate}
\end{prop}

In each panel of Figure~\ref{fig:optimal crossing point}, $F$ is globally sharper than $G$. Moreover, $0<G(z)<1$ for all $z$ in $I$, where we take $I=\RRR$ in panels (a) and (b) and $I=[0,\infty)$ in the remaining panels. So if $S$ is any of the scores listed in Proposition~\ref{prop:optimal fcst}, then the expected score $\EEE_F S(x,Y)$ exists for every $x$ in $I$, and Theorem~\ref{thm:optimal} applies. In panels (a), (c) and (d), the crossing point $\chi_A$ is unique and is therefore an optimal forecast for $S$ on $I$. However, in panels (b) and (e), there are two crossing points $\chi_A$ and $\chi_B$, each of which is a locally optimal forecast, and one of which is necessarily globally optimal.  In general, which crossing point is globally optimal depends on the measure $\mu$, as demonstrated in Table~\ref{tab:optimal depends on measure} for a selection of measures.

\begin{table}[t!]
    \renewcommand{\arraystretch}{1.5} 
    \centering \small
    \begin{tabular}{|c|c|c|}
    \hline
        $\dd\mu(\alpha)$ &  Optimal forecast in Figure~\ref{fig:optimal crossing point}b & Optimal forecast in Figure~\ref{fig:optimal crossing point}e \\ \hline
        $2\,\dd\alpha$ & $\chi_A$ & $\chi_B$ \\ 
        $\ds\frac{\dd\alpha}{\alpha(1-\alpha)}$ & $\chi_A$ & $\chi_B$  \\ 
        $\dd\delta_{1/3}(\alpha)+2\,\dd\delta_{2/3}(\alpha)$ & $\chi_B$ & $\chi_B$ \\ 
        $2\,\dd\delta_{1/3}(\alpha)+\dd\delta_{2/3}(\alpha)$ & $\chi_A$ & $\chi_A$ \\	\hline
    \end{tabular}
    \caption{The crossing point $chi_A$ or $\chi_B$ that minimizes the expected generalized diagonal score $\EEE_F\GDSd(x,Y;G,\mu)$ given the measure $\mu$ as specified, and the predictive distribution $F$ and climatological distribution $G$ of Figure~\ref{fig:optimal crossing point} (b) and (e). The measures in the first two rows yield the scoring functions $S_1$ and $S_2$, while the measures in last two rows yield specific instances of $S_5$.}
    \label{tab:optimal depends on measure}
\end{table}

Theorem~\ref{thm:optimal} provides a basis for developing a rigorous theory of the elicitability \citep[Definition~2]{gneiting2011making} of the crossing point functional, building on the work of \citet{ben2021verification}. For example, if the inequalities defining a crossing point are made strict (i.e., $F(z)<G(z)$ whenever $z<\chi$, $F(\chi)\geq G(\chi)$ and $F(z)>G(z)$ whenever $z>\chi$), then Corollary~\ref{cor:crossing point} can be strengthened to show that $\chi$ is the unique minimizer of the expected score (c.f. the end of the proof of Lemma~\ref{lem:D(x,x')}, where the corresponding inequalities become strict). A full development of this theory is beyond the scope of this paper.

\subsection{Optimal forecasts for the SEEPS score}\label{ss:SEEPS optimization}

We now discuss the implications of this theory for the SEEPS score, which admits representations in both $\GDSd$ form (Section~\ref{ss:discretize}) and $\GDSc$ form (Proposition~\ref{prop:gerrity SEEPS GDS equivalence}), depending on whether the forecast is single-valued or categorical.

When the predictive distribution has a unique crossing point, the optimal strategy is immediate: one forecasts the crossing point. However, in relevant applications (e.g., Figure~\ref{fig:optimal crossing point}e) multiple crossing points may arise. In such cases, the crossing point theorem does not determine which crossing point forecast is globally optimal.

This motivates the need for a simple and explicit decision rule. In the categorical setting, such a rule should depend only on the climatological probabilities $(p_1,p_2,p_3)$ and the predictive probabilities $(r_1,r_2,r_3)$. For single-valued forecasts, it then suffices to issue any value lying in the optimal category. A key role is played by the function $\Psi:(0,1)\times[0,1]\to\RRR$, defined by
\[
\Psi(u,v) = \frac{u-v}{u(1-u)}.
\]
The quantity $\Psi(p_i,r_i)$ measures the deviation of the predictive probability $r_i$ from the climatological probability $p_i$, normalized by the variance $p_i(1-p_i)$.

\begin{prop}\label{prop:optimal seeps}
Assume the notation above. Then optimal categorical forecasts for SEEPS are characterized as follows:
\begin{align}
C_1 \text{ is optimal } &\iff  r_1\geq p_1 \,\text{ and }\, \Psi(p_1,r_1)\leq \Psi(p_3,r_3); \\
C_2 \text{ is optimal } &\iff  r_1\leq p_1 \,\text{ and }\, r_3 \leq p_3; \\
C_3 \text{ is optimal } &\iff  r_3\geq p_3 \,\text{ and }\, \Psi(p_3,r_3)\leq \Psi(p_1,r_1).
\end{align}
\end{prop}

Thus forecasting category $C_2$ is optimal when forecast probabilities for the other categories are under their climatological probabilities. Otherwise, either $C_1$ or $C_3$ is optimal as judged by their respective forecast probabilities exceeding their climatological probabilities and having the lower normalized deviation. 

As with the generalized diagonal score, there is no universal rule that characterizes the optimal forecast in terms of crossing points when crossing points are not unique, even in the categorical setting. For example, suppose that $(p_1, p_2, p_3)=(0.7, 0.2, 0.1)$ and consider two forecasts $(r_1,r_2,r_3)=(0.75, 0.05, 0.2)$ and $(r_1',r_2',r_3')=(0.75, 0.14, 0.11)$. The differences between the cumulative predictive probabilities and the cumulative climatological probabilities yield the same sequence of signs $(+,-,0)$ for both forecasts (i.e., categories $C_1$ and $C_3$ are local crossing points), but the optimal forecast categories are $C_3$ and $C_1$ respectively.

\begin{figure}[t!]
\centering
\includegraphics[width=0.98\textwidth]{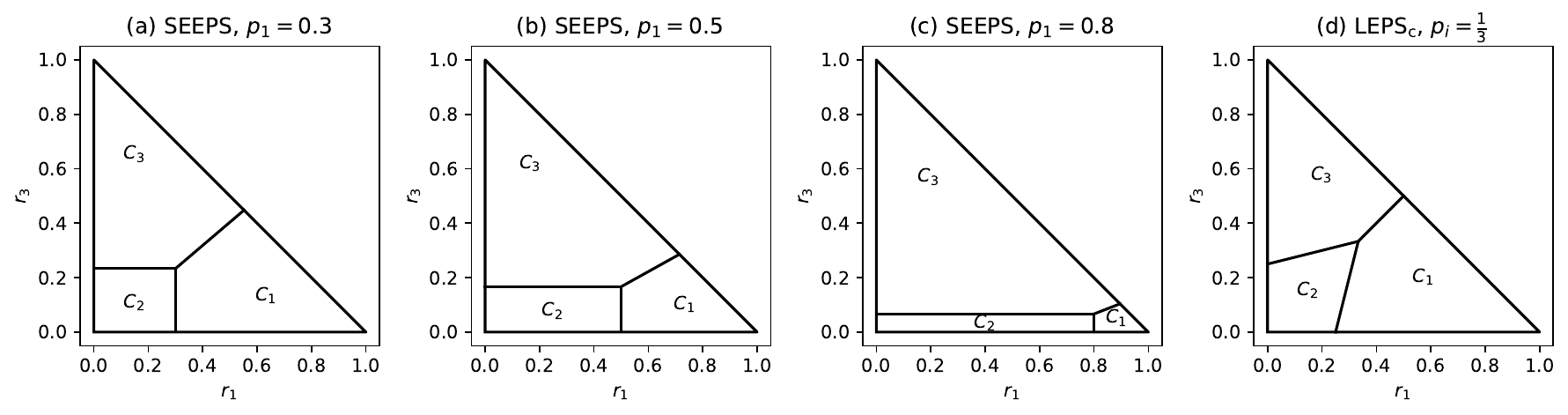}
\caption{Simplices of points $(r_1,r_3)$, where $r_i$ is the predictive probability of observing category $C_i$ for the 3-category case. Each simplex is partitioned into regions where forecasting the category $C_i$ (a-c) minimizes the expected SEEPS score for the specified climatological probability $p_1$ of dry conditions, or (d) maximizes the expected $\LEPSc$ score for equiprobable climatological probabilities ($\LEPSc$ is positively oriented). In (a-c), the boundaries of each region meet at the interior point $(p_1,p_3)$ in the simplex, and the diagonal line extending from this point has gradient $p_3(1-p_3) / (p_1(1-p_1))$. In (d), the region for $C_2$ has vertices $(0,0)$, $(1/4,0)$, $(1/3,1/3)$ and $(0,1/4)$, while the common boundary of the regions for $C_1$ and $C_3$ is the line segment from $(1/3,1/3)$ to $(1/2,1/2)$.}
\label{fig:seeps simplex}
\end{figure}

Figure~\ref{fig:seeps simplex} provides a visualization of this result when the categories are defined by \eqref{eq:SEEPS categories}, as recommended by \citet{rodwell2010new} for precipitation verification. Given predictive probabilities $r_1$ for $C_1$ (dry) and $r_3$ for $C_3$ (heavy precipitation), the point $(r_1,r_3)$ lies in one of the regions in Figure~\ref{fig:seeps simplex}a--c, where the optimal forecast category is labelled $C_i$.

For example, if $(r_1, r_3)=(0.4, 0.1)$ and the climate is relatively wet ($p_1=0.3$, Figure~\ref{fig:seeps simplex}a) then forecasting $C_1$ (dry) is optimal. By contrast, the same predictive probabilities in a middling climate ($p_1=0.5$, Figure~\ref{fig:seeps simplex}b) and a dry climate ($p_1=0.8$, Figure~\ref{fig:seeps simplex}c) yield optimal forecasts of $C_2$ (light) and $C_3$ (heavy) respectively.

\subsection{Optimal forecasts for $\LEPSc$ and the failure of crossing-point optimality}

A natural question is whether crossing-point optimality is a consequence of equitability itself, or merely of membership in the GDS class. The 3-category $\LEPSc$ score with equiprobable climatological probabilities ($p_i=1/3$), with scoring matrix (\ref{eq:LEPSCAT 3 equi}), demonstrates that crossing-point optimality does not follow from equitability alone.

Figure~\ref{fig:seeps simplex}d shows the regions where different categorical forecasts maximize the expected $\LEPSc$ score (\ref{eq:LEPSCAT 3 equi}), noting that $\LEPSc$ is positively oriented, so that larger scores are preferred. Consider predictive categorical probabilities $(r_1,r_2,r_3)=(0.05,0.65,0.3)$. For the 3-category Gerrity score with equiprobable climatological probabilities (which is in $\GDSc$ form), Corollary~\ref{cor:crossing point} implies that forecasting category $C_2$ is optimal. However, Figure~\ref{fig:seeps simplex}d shows that the optimal forecast for $\LEPSc$ is category $C_3$, which is not a crossing point. This demonstrates that equitability does not guarantee that optimal forecasts are characterized by crossing points.

We conclude with an example in which an equitable scoring rule for probabilistic forecasts is constructed from an equitable score for categorical forecasts. The Australian Bureau of Meteorology issues probabilistic seasonal forecasts for climatologically defined tercile categories. Given predictive probabilities $(r_1,r_2,r_3)$ and an observed category $i$, \citet{fawcett2005verification} assign the score
\begin{equation}\label{eq:Fawcett LEPS}
\tilde{S}((r_1,r_2,r_3),i)
= \LEPSc(1,i)r_1 + \LEPSc(2,i)r_2 + \LEPSc(3,i)r_3
\end{equation}
using the equiprobable LEPS scoring matrix (\ref{eq:LEPSCAT 3 equi}). This corresponds to the expected $\LEPSc$ score when the forecast category is selected at random according to the predictive probabilities.

It is straightforward to show that $\tilde{S}$ is equitable, a result that extends to the general case of $N$ categories and where $\LEPSc$ is replaced by any equitable score $S$ for categorical forecasts. However, the expected score under any such $\tilde{S}$ is optimized by assigning probability $1$ to the category that optimizes the expected $S$ score, and $0$ to all others.

For example, suppose that $\tilde{S}$ is given by $(\ref{eq:Fawcett LEPS})$ and that the predictive probabilities are $(0.05,0.65,0.3)$. Then category $3$ is optimal for $\LEPSc$ and consequently the probabilistic forecast $(0,0,1)$ optimizes the expected $\tilde{S}$ score.

If equitability is desired for probabilistic forecasts in categorical space, it is preferable to use proper scoring rules such as (\ref{eq:proper equitable for cats}), which do not incentivize such undesirable forecasting behaviour.

\section{Experiments}\label{s:experiments}

The preceding sections developed the generalized diagonal score framework and established its principal theoretical properties. We now explore the implications of this theory for equitability and its applications, particularly in the assessment of precipitation forecasts from global weather prediction models. To this end, we present two simple experiments.

\subsection{Experiment 1: optimal precipitation forecasts for a weather system}\label{ss:exp1}

\begin{figure}[t!]
\centering
\includegraphics[width=0.95\textwidth]{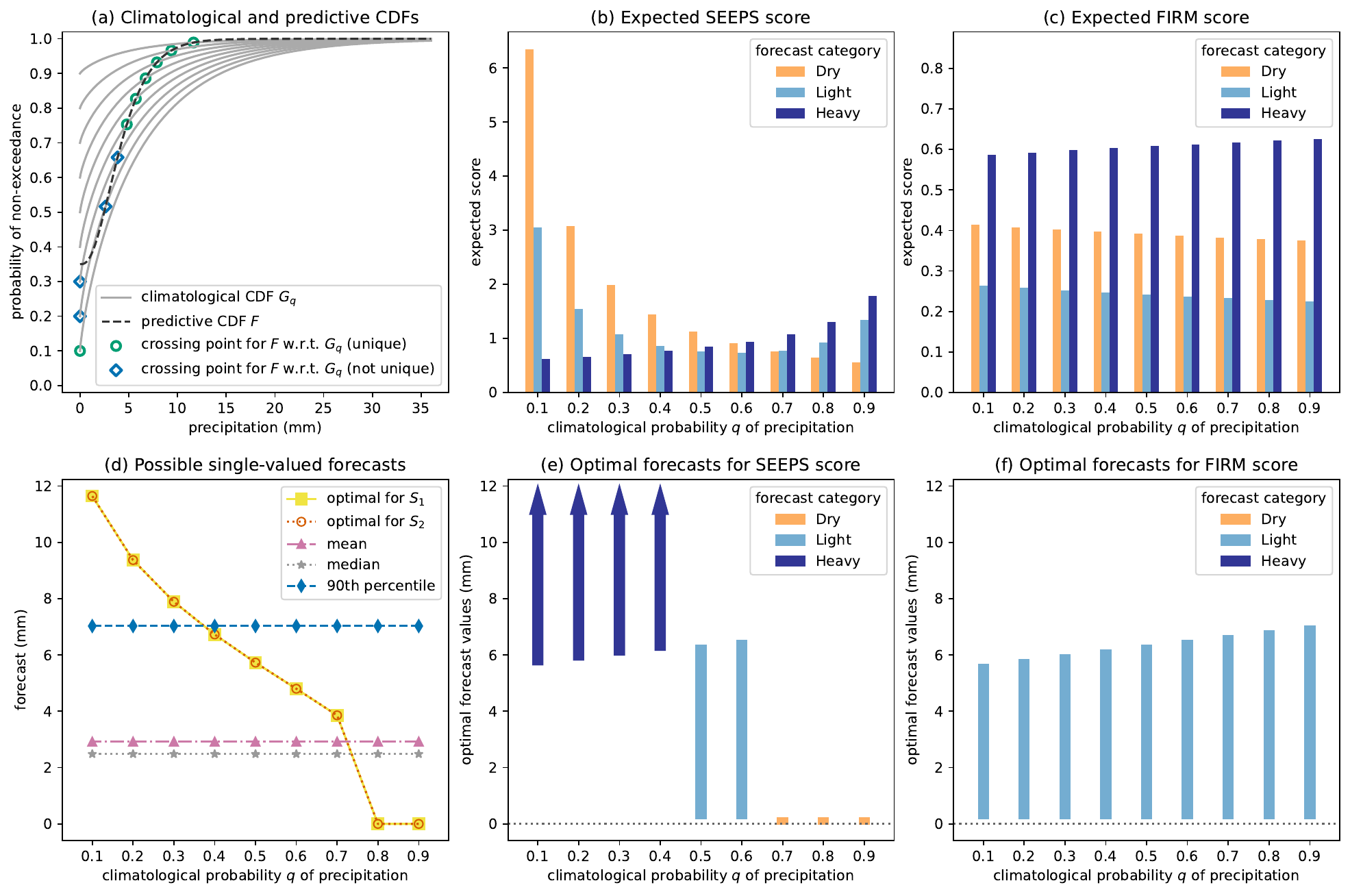}
\caption{Graphs for the thought experiment of Section~\ref{ss:exp1}. (a) Predictive CDF $F$ for precipitation (dashed dark gray curve) and climatological CDFs $G_q$ for different locations (light gray curves), along with markers representing crossing points for $F$ with respect to $G_q$; green circles when the crossing point is unique and blue diamonds otherwise. (b) Expected SEEPS score for different forecast categories, given predictive distribution $F$ and climatological probability of precipitation. A lower value is better. (c) Same as (b), but for a specific version of the FIRM score. (d) Different single-valued forecasts that a forecaster could issue given their predictive distribution $F$ and the climatological probability of precipitation at the location. (e) Range of single-valued forecasts that could be issued to optimize the forecaster's expected SEEPS score, depending on the climatological probability of precipitation at the location. (f) Same as (e), but for a specific version of the FIRM score.}
\label{fig:weather}
\end{figure}

Suppose a weather system associated with shower activity is forecast to move across a continent. At any location in its path, the forecaster assesses that the probability $F(x)$ of precipitation not exceeding $x$~mm is
\[
F(x)=0.35 + 0.65 F_\mathrm{g}(x;3,1.5) \qquad\text{whenever }x\geq0,
\]
where $F_\mathrm{g}(\,\cdot\,;\beta,\sigma)$ denotes the gamma CDF with shape $\beta$ and scale $\sigma$. The graph of $F$ is shown as the dashed black curve in Figure~\ref{fig:weather}a. For example, the probabilities of exceeding 0, 5 and 10~mm are 65\%, 23\% and 2.5\%, respectively.

The system passes over nine locations with distinct precipitation climates, ranging from very dry (climatological probability of precipitation 10\%) to very wet (90\%). The climatological CDF $G_q$ at each location is
\[
G_q(x)= (1-q) + q F_\mathrm{g}(x;0.8,2q+6) \qquad\text{whenever }x\geq0,
\]
where the probability $q$ of precipitation takes values in $\{0.1,0.2,\ldots,0.9\}$. These CDFs are shown as grey curves in Figure~\ref{fig:weather}a, with higher curves corresponding to drier climates.

We assume the forecaster is required to issue a single-valued (deterministic) forecast at each location. Possible choices include the mean (2.9~mm), median (2.5~mm), or 90th percentile (7.0~mm) of $F$, or even a random draw from $F$, which would reflect spatial variability at the cost of displacement errors. However, in the absence of a clear directive, such as a loss function to be minimized or a target functional, there is no canonical way to reduce a predictive distribution to a single-valued forecast (c.f. \citealt{murphy1985forecast, gneiting2011making}).

This lack of a canonical reduction has motivated the use of equitability as a ``desirable property'' for verification in the absence of a clear directive \citep{gandin1992equitable, Livezey2003Categorical, hogan2010equitability}. Suppose, then, that forecasts are evaluated using an equitable score derived from the $\GDS$. The theory of Section~\ref{s:optimal} implies that optimal forecasts coincide with crossing points.

Crossing points for $F$ with respect to each $G_q$ are shown in Figure~\ref{fig:weather}a: green circles indicate when the crossing point is unique, while blue diamonds indicate cases with multiple crossing points. Figure~\ref{fig:weather}d shows the resulting optimal forecasts for $S_1$ as a function of the climatological probability $q$ of precipitation, ranging from 11.6~mm in the driest climate to 0~mm in the wettest. The same forecasts are optimal under $S_2$.

These optimal forecasts vary strongly with climate despite identical predictive uncertainty across locations. This dependence arises because the optimal prediction under an equitable scoring function depends explicitly on the climatological distribution. In contrast, forecasts based on the mean, median, or other quantiles of $F$ depend only on the predictive distribution and are therefore invariant across locations (Figure~\ref{fig:weather}d).

Figure~\ref{fig:weather}b shows the expected SEEPS score for each forecast category under $F$, with categories defined by \eqref{eq:SEEPS categories}. To minimize the expected score, the forecaster should issue values in $C_3$ (`heavy precipitation') for drier climates ($q\leq 0.4$) and in $C_1$ (`dry') for wetter climates ($q\geq0.7$). Figure~\ref{fig:weather}e shows the corresponding ranges of optimal single-valued forecasts. Notably, despite identical predictive distributions across locations, the optimal forecasts depend strongly on climate: forecasts exceed roughly 6~mm in dry climates, lie between 0.2 and 6~mm for intermediate climates, and fall below 0.2~mm in wet climates. This strong dependence is primarily driven by shifts in crossing points as the climate varies, rather than by categorical boundaries adjusting with climate.

Finally, we consider a 3-category FIRM score \citep{taggart2022scoring} using the same categories (\ref{eq:SEEPS categories}) as SEEPS but with scoring matrix given by
\begin{equation}\label{eq:FIRM matrix}
\frac{1}{2}\begin{pmatrix}
0 & 1 & 2 \\
1 & 0 & 1 \\
2 & 1 & 0
\end{pmatrix}.
\end{equation}
In FIRM notation, this is $\mathrm{FIRM}((\tilde{\theta}_1,\tilde{\theta}_2), (1, 1), 0.5, 0)$, for which the optimal forecast is the category containing the median of $F$. This scoring matrix is also equivalent to the inequitable scoring matrix used as motivation by \citet[Section~3]{gandin1992equitable} to develop equitable scores.

Figures~\ref{fig:weather}c and f show the expected FIRM scores and corresponding optimal forecasts. Unlike SEEPS, the set of optimal forecasts under FIRM is relatively stable across climates, since the median of $F$ is independent of climatology. In contrast, crossing points, and hence SEEPS-optimal forecasts, vary substantially with climate.

This experiment highlights that, while the predictive distribution (and hence the mean or a quantile) is determined by meteorology (e.g., the dynamics that drive individual weather systems), optimal forecasts under equitable scoring functions are also strongly influenced by climatology. Consequently, a single-valued or categorical forecast optimized for an equitable score such as SEEPS cannot, in general, be interpreted in isolation: its interpretation depends on the local climatological distribution. For gridded forecasts, this implies that a user may require knowledge of the climatology at each grid point to interpret the forecast meaningfully. In contrast, forecasts such as the mean or median depend only on the predictive distribution and can be interpreted without reference to external climatological information.

\subsection{Experiment 2: marginal distributions of observations and forecasts}\label{ss:exp2}

A common justification for equitable scoring functions is that they encourage forecast behaviour such that the marginal distribution of the forecasts is similar to that of the observations. This experiment examines that claim under varying levels of predictive uncertainty, assuming perfectly calibrated predictive distributions.

This claim appears in several forms in the literature. For example, \citet{gandin1992equitable} argue that using certain inequitable scoring functions, such as one equivalent to the scoring matrix (\ref{eq:FIRM matrix}), for long-range monthly precipitation forecasts ``[s]uperficially ... seems quite logical and reasonable ... [but] encourages forecasters to predict the near-normal event more often than it occurs.'' The development of their equitable scoring functions is motivated, in part, by such concerns. Similar concerns motivate the equitable scoring matrices of \citet{folland1986experimental} and \citet{potts1996revised}, while \citet{rodwell2010new} explicitly claim that ``equitability also encourages `refinement' ... whereby the forecast distribution becomes equal to the observed distribution''. The experiment below examines whether such behaviour is implied by the Bayes acts of equitable scores.

Suppose that $m\in\RRR$, that $\sigma_1$ and $\sigma_2$ are positive, and that at time $t$ a random variable $Y_t$ is given by $Y_t=V_t+U_t$, where $V_t\sim \cN(m,\sigma_1^2)$, $U_t\sim \cN(0,\sigma_2^2)$ and all random variables $U_t$ and $V_t$ are pairwise independent. Then the marginal distribution of $Y_t$ (the climatological distribution) is $\cN(m, \sigma^2)$, where $\sigma^2 = \sigma_1^2 + \sigma_2^2$. Given a realization $v_t$ of $V_t$, the forecaster identifies the ideal predictive distribution $F_t = \cN(v_t, \sigma_2^2)$.

Suppose that the forecaster is required to issue both a single-valued forecast and a categorical forecast (with $\RRR$ partitioned into three categories). For single-valued forecasts, we consider two cases: (i) issuing the mean (equivalently median) of $F_t$, namely $v_t$; and (ii) issuing the crossing point $\chi_t$, which is unique and given by $\chi_t = (v_t\sigma - m\sigma_2)/(\sigma - \sigma_2)$ (see Lemma~\ref{lem:crossing point for normal}). The marginal distribution of $v_t$ is $\cN(m,\sigma_1^2)$, while the marginal distribution of $\chi_t$ is $\cN(m,(\sigma\sigma_1/(\sigma-\sigma_2))^2)$.

For categorical forecasts, we similarly consider: (i) issuing the category containing $v_t$, and (ii) issuing the category containing $\chi_t$. Forecasts of type (i) optimize the expected score under (\ref{eq:FIRM matrix}), while those of type (ii) are optimal for scoring matrices in, or equivalent to, $\GDSc$ form, including the SEEPS and Gerrity scores.

\begin{figure}[t!]
\centering
\includegraphics[width=0.95\textwidth]{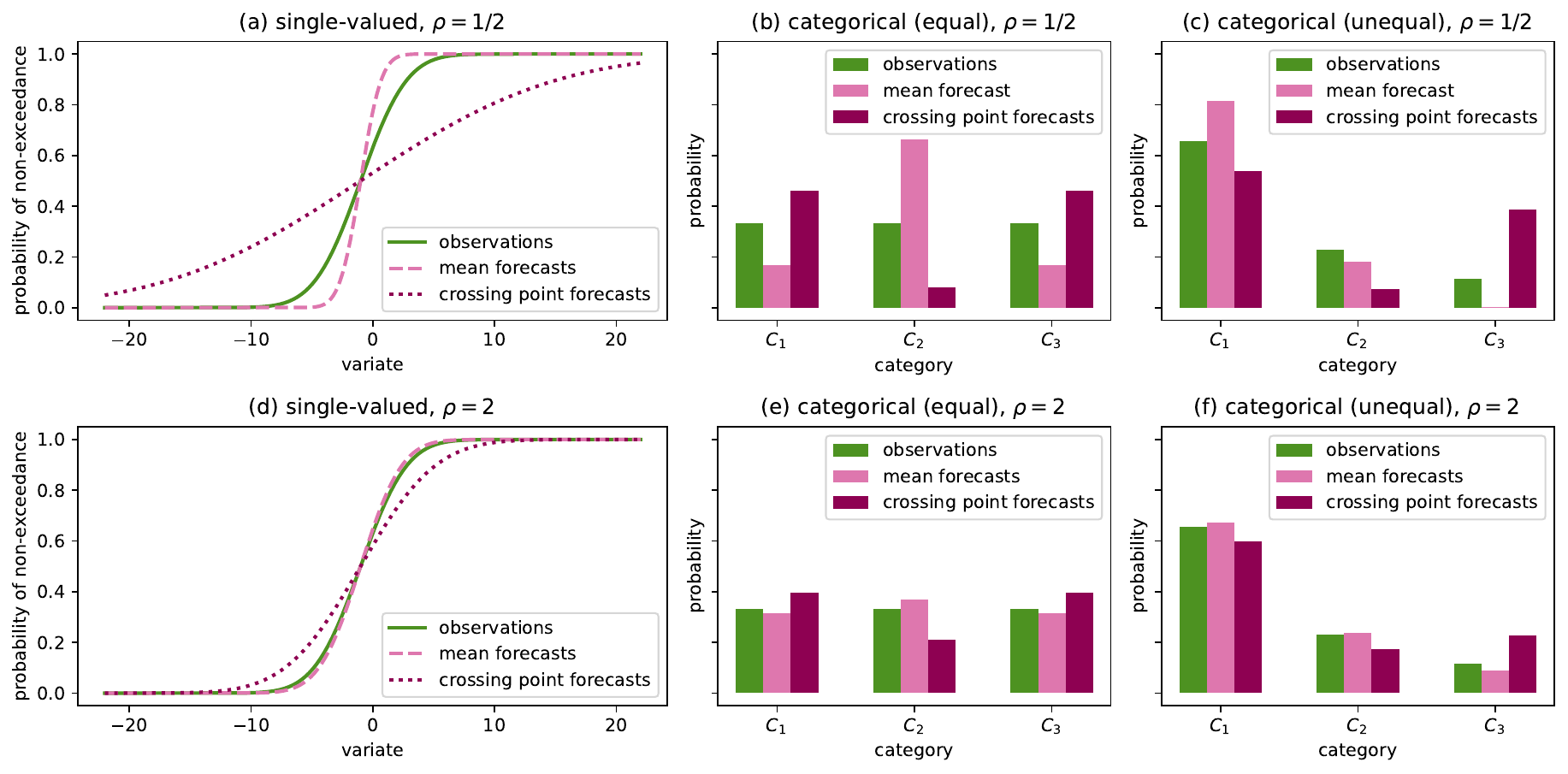}
\caption{Marginal distributions of forecasts and observations from Experiment 2, where the forecasts and observations are (a,d) single-valued, (b,e) categorical with equiprobable climatological categories, or (c,f) categorical with unequally probable climatological categories, for ideal predictive distributions with a relatively (a,b,c) low signal to noise ratio ($\rho=1/2$) or (d,e,f) high signal to noise ratio ($\rho=2$). In each case, the marginal distributions of two different types of forecasts are shown: the mean value (or category in which it lies) and the crossing point (or category in which it lies).}
\label{fig:marginal}
\end{figure}

Let $\rho$ denote the ratio $\sigma_1/\sigma_2$, which represents the signal-to-noise ratio in the predictive distribution $F_t$. Figure~\ref{fig:marginal} shows the marginal distributions of the forecasts and observations when $m=-1$, $\sigma=3$, and either $\rho=1/2$ (top row) or $\rho=2$ (bottom row).\footnote{Due to scale and translation invariance, the choice of $\mu$ and $\sigma$ is immaterial to almost all results in this section. The exception is for the 3-category case with unequal climatological probabilities, where $\mu$ and $\sigma$ were chosen to generate climatological probabilities that are somewhat typical of the categories used by SEEPS for precipitation verification as recommended by \citet{rodwell2010new}.} Forecasts and observations are treated as single-valued (first column) or converted to categories (remaining columns).

In the single-valued case (Fig.~\ref{fig:marginal}a,d), the marginal distribution of mean forecasts is narrower than that of the observations, while the distribution of crossing point forecasts is wider. For both values of $\rho$, the distribution of mean forecasts is closer to the observed marginal distribution.

In the 3-category case with equal climatological probabilities (Fig.~\ref{fig:marginal}b,e), the middle category $C_2$ is over-represented relative to observations when forecasts are based on the mean; substantially so when uncertainty is high ($\rho=1/2$), and only marginally when uncertainty is low ($\rho=2$). In contrast, crossing point forecasts under-represent $C_2$, particularly under high uncertainty.

We next consider unequal climatological probabilities (Fig.~\ref{fig:marginal}c,f) with categories defined by (\ref{eq:SEEPS categories}). Under high uncertainty, mean-based forecasts substantially under-represent $C_3$, whereas crossing-point forecasts over-represent it. Overall, mean forecasts are closer to the observed marginal distribution when uncertainty is low.

This experiment shows that forecast uncertainty is a primary driver of discrepancies between the marginal distributions of forecasts and observations. Moreover, it demonstrates that equitable scores (e.g., the Gerrity score) do not necessarily encourage forecast behaviour that brings these distributions closer together, relative to inequitable scores.

\section{Discussion}\label{s:discussion}

The paper studies optimal forecast strategies (the Bayes act) for equitable scores and clarifies their connection to crossing points for a broad class of equitable scoring functions. We now examine the implications of these results for common interpretations of equitable scores and the forecast behaviour they encourage.

\subsection{Unskillful forecasts, or forecasts with no discrimination?}\label{ss:unskillful}

At its core, a scoring function is equitable if all forecasts have the same expected score under the climatological distribution (Definition~\ref{def:equitability}). In a stationary climate, this implies that any strategy issuing constant or random forecasts achieves the same long-run mean score. However, the common paraphrase that an equitable score assigns the same expected score to all ``unskillful'' forecasts \citep{jolliffe2008impenetrable, hogan2010equitability, rodwell2010new, rodwell2011peirce} conflates skill with discrimination. We instead favour the paraphrase that an equitable score assigns the same expected score to forecasts with no discrimination. This more accurately describes characteristics of constant and random forecast strategies, and is supported by recent work \citep[Appendix~C]{bouallegue2025self}.

The distinction between skill and discrimination is important, as it recognises that calibration is a separate and essential component of forecast quality.
Recognizing this, there is no \textit{a priori} reason to require that all constant forecasts receive the same expected score if the goal is to measure predictive skill. Some constant forecasts are better calibrated than others. A more meaningful baseline of `skill' may be the expected score of the `best' constant forecast, where `best' is defined according to the application at hand, rather than treating all constant forecasts as equally uninformative.

\subsection{Equitability and bias}\label{ss:bias}

We now examine the claim that equitable scores `discourage bias in forecasting' \citep{rodwell2010new}, in contrast to inequitable scores which may `encourage forecasters---and even the developers of forecasting methods---to produce forecasts that are biased towards one event at the expense of others' \citep{gandin1992equitable}. These claims typically arise in settings where forecasts are single-valued or categorical and no directive is specified for converting predictive distributions into issued forecasts.

In such cases, if the scoring function is known, the Bayes act effectively defines the forecasting directive and thus the induced forecast behaviour. Since any reduction of a predictive distribution to a deterministic or categorical forecast necessarily introduces characteristics that may be described as `bias', it is important to clarify what type of bias is intended. As discussed in Section~\ref{ss:exp2}, differences between the marginal distributions of forecasts and observations appear to be a primary concern.

However, as demonstrated by Experiment~2, equitable scores do not uniquely mitigate such differences. For example, the Bayes act associated with the equitable Gerrity score produces different deviations in marginal distributions than the Bayes act associated with the inequitable scoring matrix (\ref{eq:FIRM matrix}). To our knowledge, there is no theoretical justification showing that equitable scores encourage, \textit{ex ante}, forecast behaviour that aligns the marginal distributions of forecasts and observations.

Experiment~1 (Section~\ref{ss:exp1}) illustrates another form of bias induced by equitable scores. For scores in GDS form (including SEEPS), optimal forecasts depend systematically on climatology. Wetter forecasts are rewarded in dry climates and drier forecasts in wet climates, even when the predictive distribution is unchanged and, in the categorical case, the thresholds defining what is `wetter' already account for climatology. In this sense, equitable scores can induce forecasts biased relative to climatology. In contrast, many inequitable scores (e.g., those consistent with the mean or a quantile) produce forecasts that depend only on the predictive distribution.

The role of bias is also relevant to the construction of the Gandin--Murphy scoring matrices, where symmetry is imposed in `view of the fact that we are concerned with the performance in a purely meteorological sense' \citep{gandin1992equitable}. Symmetry may be intended to balance penalties for over- and under-prediction. However, any symmetric matrix with unequal diagonal entries is equivalent to an asymmetric matrix with zero diagonal, obtained by subtracting each diagonal entry from its column (c.f. Proposition~\ref{prop:gerrity SEEPS GDS equivalence}(ii) and proof). Thus, symmetry alone does not ensure unbiased forecast incentives.

Ultimately, whether a scoring function encourages or discourages bias depends on the application and the intended use of the forecast. Reducing a predictive distribution to a deterministic or categorical forecast requires an implicit or explicit directive reflecting user-relevant costs and decisions. As \citet{hogan2010equitability} note, `We are then left with the problem of determining whether a directive is sensible,' after which a scoring function should be chosen that is consistent with that directive \citep{murphy1985forecast, gneiting2011making, taggart2022scoring, taggart2025warnings}. If no directive is specified, then the selected scoring function effectively defines one via its Bayes act, and its implications for forecast behaviour should be carefully considered. In particular, equitable scores should not be interpreted as rewarding forecast behaviour that eliminates bias, but rather as rewarding forecast behaviour shaped by both the predictive distribution and climatology.

\subsection{Equitability and resistance to hedging}\label{ss:hedging}

We now examine whether equitable scores discourage hedging. In the meteorological literature, hedging traditionally refers to a discrepancy between a forecaster’s `judgment' (their probabilistic beliefs, given the available information) and the issued forecast \citep{winkler1968good, murphy1978hedging}. Such hedging occurs \textit{ex ante} (before forecast issuance). Under this definition, any deterministic or categorical forecast that does not correspond to certainty under the forecaster’s judgment is, strictly speaking, hedged \citep{murphy1967note, murphy1978hedging}.

For probabilistic forecasts, hedging arises when a forecaster can improve their expected score by reporting something other than their true beliefs. Proper scoring rules preclude this behaviour: they ensure that the forecaster's expected score is optimized when they report their true probabilistic judgments. Hence, proper scoring rules actively discourage hedging. In the context of binary outcomes, \citet{jolliffe2008proper} showed that a large class of scoring functions cannot be both proper and equitable, with \citet{jolliffe2008impenetrable} concluding that `Equitability seems to be largely irrelevant for binary probabilistic forecasts'. Despite this earlier assessment, \citet{ben2018diagonal} subsequently introduced the diagonal score for predictive distributions on $\RRR$, which is both proper and equitable, while the GDS yields a broader class of scores that are proper and equitable but not strictly proper (Theorem~\ref{thm:proper not strict}), including when the outcome is binary or categorical (Corollary~\ref{cor:proper equitable for cats}). These scores do not encourage hedging, but they are insensitive to certain forms of substantial forecast misspecification, as discussed in Section~\ref{ss:optimal GDS} and illustrated in Fig.~\ref{fig:normal crossing point}.

Despite this classical definition, some authors claim that equitable scores for (non-probabilistic) categorical forecasts discourage hedging. For example, \citet{north2013assessment} describe SEEPS as `resistant to hedging' in the sense that the issued forecast differs from the one that would have been given solely to improve, \textit{ex ante}, the verification score. However, since SEEPS does not score probabilistic forecasts, this interpretation is problematic: SEEPS simply induces a reduction of the forecaster’s probabilistic judgment to a deterministic forecast via the Bayes act.

\citet{rodwell2010new} adopt a different notion of hedging, following \citet{stephenson2000use}. Here, hedging occurs if the average score over many forecast--observation pairs can be improved after the fact by systematically relabelling forecasts between categories, based on knowledge of the empirical joint distribution of forecasts and observations. This is a fundamentally different notion of hedging, constituting \textit{ex post} (after the forecast and observation are known) rather than \textit{ex ante} hedging. \citet{rodwell2010new} show that SEEPS is resistant to this form of \textit{ex post} hedging, in that relabelling a fraction of forecasts from one category $C_i$ to an adjacent category $C_j$ can improve the mean SEEPS score only when the probability of observing $C_i$, given that $C_i$ was forecast, is less than the climatological probability of $C_i$.

Nevertheless, this notion of hedging concerns only the susceptibility of a scoring function to \textit{ex post} modification of the forecasting system, and does not address the behaviour it induces \textit{ex ante}. As illustrated in Experiment~1, SEEPS induces forecast behaviour that depends strongly on climatology. In particular, it favours wetter forecasts in dry climates and drier forecasts in wet climates, even when the underlying predictive information is the same. In effect, SEEPS rewards a model that captures atmospheric dynamics well but exhibits a precipitation scheme that is relatively overactive in dry climates and underactive in wet climates.

\subsection{Scoring precipitation forecasts from global weather prediction models}\label{ss:alternatives}

Finally, we consider alternative approaches to evaluating single-valued precipitation forecasts from global weather prediction models. Given the diversity of forecast users and applications, in the face of uncertainty no single-valued precipitation forecast will satisfy all requirements. Nonetheless, each candidate scoring function rewards particular forecast behaviours, so it is important to understand what these are and whether those behaviours are appropriate for forecasts that may be used across a wide range of applications. One approach is to identify a statistical functional with desirable properties in this context, and then consider which consistent scoring functions  align with additional evaluation objectives.

To illustrate this approach, we consider the expectation, which has several advantages over some alternative functionals for precipitation forecasting. First, expectations are additive across time and space: accumulated precipitation over multiple time steps or across a catchment is itself an expectation, regardless of the complexity of the spatio-temporal dependence structure. In contrast, analogous aggregations of medians or crossing points do not admit a simple interpretation. Second, the expectation reflects changes in both the likelihood and intensity of precipitation, whereas some functionals, such as the median, can remain unchanged over a wide range of circumstances. As a result, expectation forecasts may discriminate between situations that are treated identically by some other summaries, making them potentially more useful inputs to post-processing systems. Third, the expectation for a particular time and location depends only on the distribution of possible precipitation amounts, given information available to the model, whereas the forecasts preferred by equitable scores depend on both this distribution and local climatology. Consequently, spatial variation in expected precipitation reflects variation in the forecast weather systems and their interaction with orography, whereas forecasts optimized for equitable scores may also reflect regional differences in precipitation climatology.

Scoring functions $S$ that are consistent for the expectation functional are precisely those of \textit{Bregman form} (\citealt{savage1971elicitation}; \citealt[Theorem~7]{gneiting2011making}):
\begin{equation}\label{eq:bregman}
S(x,y) = \phi(y) - \phi(x) - \phi'(x)(y-x),
\end{equation}
where $\phi$ is a convex function with subderivative $\phi'$ (the subderivative is a generalisation of the derivative, and equals the derivative when $\phi$ is differentiable). If $\phi$ is strictly convex then $S$ is strictly consistent for the expectation. A well-known example of a strictly consistent scoring function is the squared error $S(x,y)=(x-y)^2$, obtained by taking $\phi(t)=t^2$. However, some have criticized this scoring function for being overly sensitive to large errors, causing verification statistics to be dominated by extreme events. An alternative is to use the climatological distribution $G$ to construct a strictly consistent scoring function that is less sensitive to large errors. For example, when $G$ is a suitably chosen parametric distribution fitted to observational data, such as $G_2$ of Figure~\ref{fig:real valued equitable scores}, then $\phi(t)=-G(t)$ yields the strictly consistent scoring function
\begin{equation}\label{eq:bregman using climate}
S(x,y)=G(x)-G(y)+G'(x)(y-x), \qquad\forall x,y\in[0,\infty).
\end{equation}
This is less sensitive to large errors, since $G'(x)\to 0$ as $x\to\infty$ and $|G(x)-G(y)|\leq 1$. A similar scoring function could be constructed using the climatological distribution conditioned on precipitation occurrence. Additionally, scaling the score by a suitable transformation of, for example, mean seasonal precipitation at each location can reduce the influence of predictive performance in high-precipitation regions on aggregated global verification results.

Such strictly consistent scores reward increasingly accurate estimates of the expectation and may therefore favour relatively smooth forecast fields, potentially benefiting data-driven models relative to physics-based models \citep{pan2025comparative, sha2025investigating}. If this is undesirable, one can consider using scores that are consistent but not strictly consistent. For example, the family of FIRM scores \citep{taggart2022scoring} contains scoring functions that are consistent for the expectation while retaining a categorical structure, similar to SEEPS. Given categorical thresholds $\theta_1$ and $\theta_2$ as in (\ref{eq:SEEPS categories}), and positive weights $w_1$ and $w_2$, consider the scoring function
\begin{equation}\label{eq:FIRM expectation}
S(x,y) = w_1 S^\mathrm{E}_{\theta_1}(x,y) + w_2 S^\mathrm{E}_{\theta_2}(x,y),
\end{equation}
where
\[
 S^\mathrm{E}_\theta(x,y)=
\begin{cases}
\tfrac12 (\theta-y), & y\leq \theta < x, \\
\tfrac12 (y-\theta), & x\leq \theta < y, \\
0, & \text{otherwise}.
\end{cases}
\]
In FIRM notation, $S$ corresponds to $\mathrm{FIRM}((\theta_1,\theta_2),(w_1,w_2),0.5,\infty)$, for which any forecast lying in the category that contains the expected value is optimal. In particular, $S$ is of Bregman form \eqref{eq:bregman} via the convex function $\phi(t) = (w_1\max(t,\theta_1)+w_2\max(t,\theta_2))/2$. Were this score $S$ used for Experiment 1, the corresponding ranges of optimal forecasts match those illustrated in Figure~\ref{fig:weather}f. The weights $w_i$ reflect the relative importance of discrimination across thresholds. For example, choosing $w_2=10w_1$ prioritizes correct classification relative to the heavy precipitation threshold $\theta_2$ over the dry threshold $\theta_1$. Evaluation using FIRM can be readily compared with SEEPS using the implementations in the open-source \texttt{scores} package \citep{leeuwenburg2024scores}.

The expectation provides only one example of this general approach. However, if median or $\alpha$-quantile forecasts are considered more appropriate for a particular application, an analogous theory of consistent scoring functions exists for those functionals \citep[Theorem~3.2]{gneiting2011quantiles}, including the $\mathrm{FIRM}((\theta_1,\theta_2),(w_1,w_2),\alpha,0)$ score for discretization into three categories \citep{taggart2022scoring}.  Unlike the corresponding scores for the expectation, some of these analogous scores are bounded on $[0,\infty)\times[0,\infty)$. Examples for the median include $S(x,y)=|G(x)-G(y)|$ and the FIRM scoring matrix \eqref{eq:FIRM matrix} with weights $w_1=w_2=1$.

\subsection{Summary and recommendations}\label{ss:summary}

The discussion in this section can be summarized in the following recommendations.
\begin{enumerate}
\item To describe equitability more precisely, we advocate avoiding the term `unskillful' in this context. Instead, we recommend the formulation that `an equitable score assigns the same expected score to all forecast systems with no discrimination.' This recognises that calibration is a distinct component of predictive skill. Further work on the relationship between equitability and discrimination, building on \citet[Appendix~C]{bouallegue2025self}, would be valuable.
\item We advocate abandoning the view that equitable scores for single-valued or categorical forecasts discourage biased forecasts, unless both a clearly defined notion of bias and compelling arguments supporting this claim are provided. A more informative description is that many equitable scores, including widely used ones such as SEEPS and the Gerrity score, encourage forecasts that coincide with crossing points (Definition~\ref{def:crossing point}), where the predictive distribution transitions from below the climatological distribution to above it.
\item Given the increasingly varied use of the term `hedging' in the meteorological literature and beyond \citep{jolliffe2008proper}, we recommend instead using terms with well-established meanings, such as `proper scoring rule' and `consistent scoring function' \citep{murphy1985forecast, gneiting2007strictly, gneiting2011making}. While the generalized diagonal score for probabilistic forecasts is both proper and equitable, we discourage describing equitable scores for categorical or single-valued forecasts as resistant to hedging. Instead, it is more informative to describe explicitly how score optimization depends on the predictive and climatological distributions, for example through crossing points.
\item When verifying single-valued precipitation forecasts from global weather prediction models, practitioners should carefully consider the type of forecast behaviour they wish to encourage and select scoring functions accordingly. In particular, the suitability of SEEPS depends on whether the climatology-dependent forecast behaviour it rewards, as illustrated in Section~\ref{ss:exp1}, is desirable for the intended application.
\end{enumerate}
These recommendations aim to clarify the interpretation of equitability and to support more informed choices about the use of equitable scores in practice.

\vspace{0.2cm}

\noindent\textbf{Acknowledgements.} The authors thank Deryn Griffiths and Mohammadreza Khanarmuei for reviewing an earlier version of this manuscript.

\vspace{0.2cm}

\noindent\textbf{Data availability statement.} A repository containing the data and code used for this paper is available online \citep{taggart2026code}.

\vspace{0.2cm}

\noindent\textbf{Conflict of interest statement.} The authors declare no conflict of interest.

\appendix

\section{Proofs of results}\label{a:proofs}

\begin{lem}\label{lem:alpha in RG}
Suppose that $G\in\mathcal{P}$ and $\alpha\in\RG$
\begin{enumerate}
\item[(i)] Then $G(G^{-1}(\alpha))=\alpha$.
\item[(ii)] Suppose also that $G$ is strictly increasing on $\mathcal{D}(G)$. Then $G^{-1}(\alpha) < x$ if and only if $\alpha < G(x)$.
\end{enumerate}
\end{lem}

\begin{proof}
To prove (i), note the standard fact that $G(G^{-1}(\alpha))\geq \alpha$ whenever $0<\alpha<1$, which follows from the right-continuity of $G$ and the definition of its inverse. On the other hand, if $\alpha\in\RG$ then $\alpha=G(z_0)$ for some $z_0$ and hence $z_0$ belongs to the set $\{z\in\RRR: G(z)\geq \alpha\}$. It follows that $G^{-1}(\alpha)=\min\{z\in\RRR: G(z)\geq \alpha\}\leq z_0$. Applying $G$ to both sides and using the fact that $G$ is nondecreasing gives $G(G^{-1}(\alpha))\leq G(z_0)=\alpha$.

To prove (ii), suppose that $G$ is strictly increasing. If $G^{-1}(\alpha) < x$ then applying $G$ to both sides and using part (i), $\alpha = G(G^{-1}(\alpha)) < G(x)$, which establishes the forward implication. To prove the converse, suppose that $\alpha<G(x)$. Then $\alpha\leq G(x)$ whence
\begin{equation}\label{eq:RG lemma}
G^{-1}(\alpha) \leq x
\end{equation}
by the definition of the generalized inverse of $G$. Now if equality holds in \eqref{eq:RG lemma} then $\alpha=G(G^{-1}(\alpha))=G(x)$, which contradicts the premise that $\alpha<G(x)$. Hence $G^{-1}(\alpha) < x$, completing the proof.
\end{proof}

\begin{proof}[Proof of Proposition~\ref{prop:elementary score equitable}]
Suppose that $G\in\mathcal{P}$, $\alpha\in(0,1)$ and $\theta\in I$. Note that
\begin{align}
\EEE_G S_{\alpha,\theta}(x,Y) 
&= \alpha\PPP_G(x\leq \theta<Y) + (1-\alpha)\PPP_G(x> \theta\geq Y) \notag \\
&= \alpha(1-G(\theta))\one\{x \leq \theta\} + (1-\alpha)G(\theta)\one\{x > \theta\}. \label{eq:exp ES}
\end{align}

Now suppose that $S_{\alpha,\theta}$ is equitable with respect to $G$. Then there exists a constant $c$ such that the right hand side of (\ref{eq:exp ES}) equals $c$ for every $x$ in $I$. This implies that $\alpha(1-G(\theta)) =  (1-\alpha)G(\theta)$, which in turn implies that $G(\theta)=\alpha$.

To prove the converse, suppose that $G(\theta)=\alpha$.  Then (\ref{eq:exp ES}) implies that $\EEE_G S_{\alpha,\theta}(x,Y)=\alpha(1-\alpha)$ for every $x$ in $I$.
\end{proof}

\begin{proof}[Proof of Corollary~\ref{cor:diagonal score equitable}]
Suppose that $S_{\alpha,G^{-1}(\alpha)}$ is equitable. Then the forward direction of  Proposition~\ref{prop:elementary score equitable} gives $G(G^{-1}(\alpha))=\alpha$, whence $\alpha\in\RG$. To prove the converse, suppose that $\alpha\in\RG$. Lemma~\ref{lem:alpha in RG} (i) implies that $G(G^{-1}(\alpha))=\alpha$ and the backward direction of Proposition~\ref{prop:elementary score equitable} implies that  $S_{\alpha,G^{-1}(\alpha)}$ is equitable.
\end{proof}

\begin{proof}[Proof of Theorem~\ref{thm:GDS is equitable}]
Suppose that $F\in\mathcal{P}_0$. Then
\begin{align*}
\EEE_G\GDS(F,Y;G,\mu)
&= \int_{\RG} \EEE_G S_{\alpha,G^{-1}(\alpha)}(F^{-1}(\alpha),Y) \, \dd\mu(\alpha) \\
&= \int_{\RG} \alpha(1-\alpha) \, \dd\mu(\alpha),
\end{align*}
where the first equality is justified by Fubini's theorem and the second by Proposition~\ref{prop:elementary score equitable}.
\end{proof}

\begin{lem}\label{lem:finite expectations}
If $G\in\mathcal{P}$ and $0<G(x)<1$ for all $x$ in $\RRR$ then $0\leq \EEE_G[-\ln(G(Y))]\leq 1$ and $0\leq\EEE_G[-\ln(1-G(Y))]\leq 1$.
\end{lem}

\begin{proof}
It is clear that $-\ln(G(Y))> 0$ almost surely and hence the expectation has a lower bound of $0$. To prove the upper bound, the tail-probability expectation formula gives
\begin{align*}
\EEE_G[-\ln(G(Y))]
&= \int_0^\infty \PPP_G(-\ln(G(Y)) > t) \,\dd t \\
&= \int_0^\infty \PPP_G(G(Y) < \mathrm{e}^{-t}) \,\dd t \\
&\leq \int_0^\infty \mathrm{e}^{-t} \,\dd t \\
&=1,
\end{align*}
where the inequality follows from the standard fact of the probability integral transform that $\PPP(G(Y)\leq u)\leq u$ whenever $0<u\leq 1$. This proves the result for $\EEE_G[-\ln(G(Y))]$. The proof for $\EEE_G[-\ln(1-G(Y))]$ is similar.
\end{proof}

\begin{proof}[Proof of Proposition~\ref{prop:real valued equitable scores}]
To prove (i), suppose that $G$ satisfies either assumption (G1) or (G2) and the measure $\mu$ has no atoms. The continuity assumptions give $\RG=(0,1)$ under (G1) and $\RG=[G(0),1)$ under (G2). The strict monotonicity assumption on $G$ and Lemma~\ref{lem:alpha in RG} (i) gives
\begin{align*}
\{\alpha\in\RG:x\leq G^{-1}(\alpha)<y\}
&=\{\alpha\in\RG:G(x)\leq \alpha< G(y)\} \\
&=[G(x),G(y))
\end{align*}
whenever $x,y \in I$ and $x<y$. Then \eqref{eq:deterministic GDS} gives
\[
\GDSd(x,y;G,\mu)
=
\begin{cases}
\displaystyle\int_{[G(x),G(y))}  \alpha \,\dd\mu(\alpha), & x<y, \\
\displaystyle\int_{[G(y),G(x))} (1-\alpha) \,\dd\mu(\alpha), & x>y, \\
0, & x=y.
\end{cases}
\]
Equation~\eqref{eq:GDSd reduced} follows from the hypothesis that $\mu$ has no atoms.

To prove (ii), let $\dd\mu(\alpha)=2\,\dd\alpha$ and suppose that $G$ satisfies (G1) or (G2). Then by part (i),
\[
\GDSd(x,y;G,\mu)= \int_{G(y)}^{G(x)}2(\one\{x>y\} - \alpha)\,\dd\alpha = S_1(x,y).
\]
To prove equitability, note that $S_1(x,y)$ is bounded by $1$ and consequently the expected score $\EEE_G S_1(x,Y)$ is also bounded, so Theorem~\ref{thm:GDS is equitable} applies, and the equitability constant is readily calculated from \eqref{eq:GDS equitability constant} using $\RG=(0,1)$ under assumption (G1) or  $\RG=[G(0),1)$ under (G2). The equivalence of $S_1$ and $\LEPSr$ follows from the identity
\begin{equation}\label{eq:S1 LEPS equivalence}
S_1(x,y)=-\tfrac{1}{3}\,\LEPSr(x,y) + 2G(y)^2-2G(y)+\tfrac{2}{3}.
\end{equation}

To prove (iii), let $\dd\mu(\alpha)=(\alpha(1-\alpha))^{-1}\,\dd\alpha$ and suppose that $G$ satisfies (G3) or (G4). Then by part (i),
\[
\GDSd(x,y;G,\mu) = \int_{G(y)}^{G(x)} \frac{\dd\alpha}{\one\{x>y\} - \alpha} = S_2(x,y).
\]
To establish equitability, fix $x$ in $\RRR$ and note that
\begin{equation}\label{eq:S2 decomp}
\EEE_G S_2(x,Y)=\EEE_G \ln|\one\{x< Y\} - G(x)| + \EEE_G[-\one\{x< Y\}\ln(1-G(Y))] + \EEE_G[-\one\{x \geq Y\}\ln(G(Y))].
\end{equation}
By Lemma~\ref{lem:finite expectations} each term on the right hand side is finite, and hence by Theorem~\ref{thm:GDS is equitable} $S_2$ is equitable with respect to $G$. The equitability constant follows from \eqref{eq:GDS equitability constant} with $\RG=(0,1)$ under (G3) and $\RG=[G(0),1)$ under (G4). The identity $S_2(x,y)= -\Grin(x,y) - \ln(G(y)) - \ln(1-G(y)) - 1$ establishes that $S_2$ and $\Grin$ are equivalent.

The proof of (iv) proceeds similarly to (ii), noting that $\mu$ has no atoms and $\EEE_G S_3(x,Y)<\infty$.
\end{proof}

\begin{proof}[Proof of Proposition~\ref{prop:GDSd ecdf}]
Note that $\RG=\{P_i\}_{i=1}^{N-1}$, $\mu$ is concentrated on $\RG$, $G^{-1}(P_i)=z_i$ and the integrand of \eqref{eq:deterministic GDS} is precisely $S_{\alpha,G^{-1}(\alpha)}(x,y)$. The expected score with respect to $G$ is clearly finite, so Theorem~\ref{thm:GDS is equitable} applies.
\end{proof}

\begin{proof}[Proof of Corollary~\ref{cor:GDSc}]
This is immediate from Proposition~\ref{prop:GDSd ecdf} and the construction of $\GDSc$.
\end{proof}

\begin{proof}[Proof of Proposition~\ref{prop:gerrity SEEPS GDS equivalence}]
The proof of part (i) is straightforward. To prove part (ii), let $g_{i,j}$ and $h_{i,j}$ denote the $(i,j)$th entry of the Gerrity scoring matrix (\ref{eq:Gerrity 3-cat}) and SEEPS scoring matrix (\ref{eq:SEEPS}) respectively. Then a simple calculation shows that $h_{i,j}=g_{j,j}-g_{i,j}$ for every $i$ and $j$. 

We now turn to the proof of part (iii). The scoring matrix $(g_{i,j})_{i,j=1}^N$ for the Gerrity score is defined by
\begin{align*}
g_{i,i}&=\frac{1}{N-1}\Big(\sum_{k=1}^{i-1}R(k)+ \sum_{k=i}^{N-1}D(k)\Big), & 1\leq  i \leq N, \\
g_{i,j}&=\frac{1}{N-1}\Big(\sum_{k=1}^{i-1}R(k)+ \sum_{k=i}^{j-1}(-1) + \sum_{k=i}^{N-1}D(k)\Big), & 1\leq i < j \leq N, \\
g_{i,j}&=g_{j,i}, & 1\leq j < i \leq N, \\
\intertext{where}
D(k) &= \Big(1-\sum_{r=1}^k p_r \Big) / \sum_{r=1}^k p_r, & 1\leq  k \leq N-1, \\
R(k) &= 1/D(k), & 1\leq  k \leq N-1
\end{align*}
\citep{gerrity1992note}. Let $m_{i,j}$ denote the $(i,j)$ entry of the scoring matrix for $\GDSc(\,\cdot\,,\,\cdot\,; \boldP, \boldw)$ as per the hypothesis of the theorem, with weights given by $w_k=((N-1)P_k(1-P_k))^{-1}$. It will be shown that $m_{i,j}=g_{j,j}-g_{i,j}$, which implies the result.

First, if $1\leq i < j \leq N$ then
\begin{align*}
g_{j,j}-g_{i,j}
&=\frac{1}{N-1}\sum_{k=i}^{j-1}(R(k)+1) \\
&= \frac{1}{N-1}\sum_{k=i}^{j-1}\left(\frac{P_k}{1-P_k}+1\right) \\
&= \sum_{k=i}^{j-1}\frac{1}{(N-1)(1-P_k)} \\
&= \sum_{k=i}^{j-1}w_k P_k \\
&= m_{i,j}.
\end{align*}
Second, if $1\leq j < i \leq N$ then
\begin{align*}
g_{j,j}-g_{i,j}
&=\frac{1}{N-1}\sum_{k=j}^{i-1}(D(k)+1)\\
&= \frac{1}{N-1}\sum_{k=j}^{i-1}\left(\frac{1-P_k}{P_k}+1\right) \\
&= \sum_{k=j}^{i-1} \frac{1}{(N-1)P_k} \\
&= \sum_{k=j}^{i-1} w_k (1-P_k) \\
&= m_{i,j}.
\end{align*}
Third, if $1\leq  i \leq N$ then $g_{i,i}-g_{i,i}=0=m_{i,i}$, which completes the proof.
\end{proof}

\begin{proof}[Proof of Proposition~\ref{prop:LEPSCAT}]
Let $\mathcal{U}(a,b)$ denote the uniform distribution on an interval with endpoints $a$ and $b$. Since $G$ is continuous, the probability integral transform implies that (\ref{eq:LEPSCAT}) equals
\begin{equation}\label{eq:LEPSCAT PIT}
\EEE[3(1-|U-V| + U^2 - U + V^2 - V) - 1  | U\sim \mathcal{U}(P_{i-1},P_i), V\sim \mathcal{U}(P_{j-1},P_j)],
\end{equation}
where $U$ and $V$ are independent. Suppose that $U\sim \mathcal{U}(P_{i-1},P_i)$ and $V\sim \mathcal{U}(P_{j-1},P_j)$ are independent. Then $\EEE U=(P_{i-1}+P_i)/2$ and  $\EEE U^2=(P_{i-1}^2+P_{i-1} P_i+P_i^2)/3$. Expectations for terms involving $V$ only can be obtained from these by exchanging $U$ for $V$ and $i$ for $j$. 
If $i>j$ then $U \geq V$ almost surely, and hence $\EEE|U-V|=\EEE U - \EEE V = (P_{i-1}+P_i)/2 - (P_{j-1}+P_j)/2$. If $i=j$ then $\EEE|U-V|=(P_i-P_{i-1})/3$. 
Substituting these calculations into (\ref{eq:LEPSCAT PIT}) and simplifying gives (\ref{eq:LEPSCAT entry 1}) and (\ref{eq:LEPSCAT entry 2}). Equation~(\ref{eq:LEPSCAT entry 3}) follows from symmetry.

To prove equitability, let $Y_\mathrm{c}$ denote the random variable taking values in categorical space $\{C_1,\ldots,C_N\}$, and $G_\mathrm{c}$ its climatological distribution, such that $\PPP_{G_\mathrm{c}}(Y_\mathrm{c}=C_i)=p_i$. Let $G_i$ denote the conditional distribution of the real-valued random variable $Y$ given $Y\in C_i$. Let $c$ denote the equitability constant of $\LEPSr$. Then whenever $1\leq i\leq N$,
\begin{align*}
&\EEE [\LEPSc(C_i, Y_\mathrm{c}) \mid Y_\mathrm{c} \sim G_\mathrm{c}] \\
&\quad=\sum_{j=1}^N p_j\,\LEPSc(C_i, C_j) \\
&\quad=\sum_{j=1}^N p_j\,\EEE[\LEPSr(X,Y) \mid X\sim G, Y\sim G, X\in C_i, Y\in C_j] \\
&\quad=\int_{C_i} \Big(\sum_{j=1}^N p_j \int_{C_j} \LEPSr(x,y)\,\dd G_j(y) \Big) \,\dd G_i(x) \\
&\quad=\int_{C_i} \Big(\int_I \LEPSr(x,y)\,\dd G(y) \Big) \,\dd G_i(x) \\
&\quad=\int_{C_i} \EEE_G\LEPSr(x,Y) \,\dd G_i(x) \\
&\quad=\int_{C_i} c \,\dd G_i(x) \\
&\quad=c.
\end{align*}
Finally, $c=0$ by the construction of $\LEPSr$ \citep{potts1996revised}, or equivalently by taking the expectation of $Y$ with respect to $G$ in (\ref{eq:S1 LEPS equivalence}) and using the equitability constant $1/3$ of $S_1$.
\end{proof}

\begin{proof}[Proof of Proposition~\ref{prop:gandin murphy not GDS}] Write $Q \equiv Q'$ if two scoring matrices $Q$ and $Q'$ are equivalent, and denote the matrix (\ref{eq:LEPSCAT 3 equi}) by $Q$. Then 
\begin{align}
Q 
&\equiv -9Q \notag \\
&\equiv
\begin{pmatrix}
-8 & 1 & 7 \\
1 & -2 & 1 \\
7 & 1 & -8
\end{pmatrix}
+
\begin{pmatrix}
8 & 2 & 8 \\
8 & 2 & 8 \\
8 & 2 & 8
\end{pmatrix} \label{eq:add matrix} \\
&=
\begin{pmatrix}
0 & 3 & 15 \\
9 & 0 & 9 \\
15 & 3 & 0
\end{pmatrix}, \label{eq:final matrix}
\end{align}
where the last matrix in (\ref{eq:add matrix}) is the unique matrix whose entries depend only on the observed category (i.e., columns) so that $-9Q$ is transformed into an equivalent matrix with zeros on the diagonal and positive entries otherwise. Suppose now that the matrix of (\ref{eq:final matrix}) is equivalent to the matrix of (\ref{eq:GDS scoring matrix 3-cat}). Then, by matching entries, there must be positive numbers $w_1$, $w_2$, $P_1$ and $P_2$ such that $w_1 P_1=3$, $w_2 P_2=9$ and $w_1 P_1+w_2 P_2=15$. This is impossible, so the assumption that $Q$ is equivalent to the matrix of (\ref{eq:GDS scoring matrix 3-cat}) is false. The same method of proof can be used for the other scoring matrices listed in the proposition.
\end{proof}

\begin{proof}[Proof of Proposition~\ref{prop:equitable discretize}]
Assume the hypotheses of the proposition. Then $G(\theta_i)\in\RG$ for every $i=1,\ldots,N-1$. Let $\mu=\sum_{i=1}^{N-1} w_i \delta_{G(\theta_i)}$. Then the scoring function
\[
\GDSd(x,y;G,\mu)=\sum_{i=1}^{N-1} w_i S_{G(\theta_i), G^{-1}(G(\theta_i))}(x,y).
\]
is equitable by Theorem~\ref{thm:GDS is equitable}, since its expected score given $Y\sim G$ is bounded by $\sum_{i=1}^{N-1}w_i$. The hypothesis $G^{-1}(G(\theta_i)) = \theta_i$ implies that $\GDSd(x,y;G,\mu)=S_5(x,y)$. Finally, $S_5(x,y)$, being a weighted sum of elementary scores, is constant on each set $C_j \times C_k$ and can change value only when $x$ or $y$ crosses one of the thresholds $\theta_i$. Hence $S_5(x,y)$ depends only on the categories to which $x$ and $y$ belong.
\end{proof}

\begin{proof}[Proof of Theorem~\ref{thm:proper not strict}]
The scoring function $S_{\alpha,G^{-1}(\alpha)}:\RRR\times\RRR\to\RRR$ is consistent for the $\alpha$-quantile functional relative to the class $\mathcal{P}_0$ \citep{gneiting2011quantiles}, whence the scoring rule $F\times y\mapsto S_{\alpha,G^{-1}(\alpha)}(F^{-1}(\alpha),y)$ is proper relative to the class $\mathcal{P}_0$ \citep[Theorem~3]{gneiting2011making}. Hence, if $F,\tilde{F}\in\mathcal{P}_0$ then by Fubini's theorem,
\begin{align}
\EEE_F \GDS(\tilde{F},Y;G,\mu)
&= \int_{\RG} \EEE_F S_{\alpha,G^{-1}(\alpha)}(\tilde{F}^{-1}(\alpha),Y)\,\dd\mu(\alpha) \notag\\
& \geq  \int_{\RG} \EEE_F S_{\alpha,G^{-1}(\alpha)}(F^{-1}(\alpha),Y)\,\dd\mu(\alpha) \notag\\
&= \EEE_F \GDS(F,Y;G,\mu). \label{eq:proper proof}
\end{align}
This proves part (i). Part (ii) follows immediately from the equitability of $\GDS(\twocdots;G,\mu)$ with respect to $G$. To prove part (iii), suppose that (\ref{eq:same sign condition}) holds. Then whenever $\alpha\in\RG$,
\begin{align*}
&\EEE_F S_{\alpha,G^{-1}(\alpha)}(\tilde{F}^{-1}(\alpha),Y) \\
&\quad= \alpha \PPP_F(\tilde{F}^{-1}(\alpha) \leq G^{-1}(\alpha) < Y) + (1-\alpha) \PPP_F(\tilde{F}^{-1}(\alpha) > G^{-1}(\alpha) \geq Y) \\
&\quad= \alpha\one\{\tilde{F}^{-1}(\alpha) \leq G^{-1}(\alpha)\}(1-F(G^{-1}(\alpha))) +(1-\alpha) \one\{\tilde{F}^{-1}(\alpha) > G^{-1}(\alpha) \} F(G^{-1}(\alpha)) \\
&\quad= \alpha\one\{F^{-1}(\alpha) \leq G^{-1}(\alpha)\}(1-F(G^{-1}(\alpha))) +(1-\alpha) \one\{F^{-1}(\alpha) > G^{-1}(\alpha) \} F(G^{-1}(\alpha)) \\
&\quad=  \EEE_F S_{\alpha,G^{-1}(\alpha)}(F^{-1}(\alpha),Y),
\end{align*}
whereby one obtains equality in (\ref{eq:proper proof}).
\end{proof}

\begin{lem}\label{lem:F on integers}
Suppose that $F$ is a distribution in $\mathcal{P}$ concentrated on the integers, $i$ is an integer and $0<\alpha<1$. Then $i<F^{-1}(\alpha)$ if and only if $F(i)<\alpha$.
\end{lem}

\begin{proof}
The forward implication is a direct consequence of the definition of $F^{-1}$. Conversely, suppose that $F(i)<\alpha$. Since $F$ is concentrated on the integers and right-continuous, the smallest possible $x$ in $\RRR$ such that $F(x)\geq \alpha$ is $i+1$. That is, $F^{-1}(\alpha)\geq i+1>i$.
\end{proof}

\begin{proof}[Proof of Corollary~\ref{cor:proper equitable for cats}]
Suppose that $F,G\in\mathcal{P}_\mathrm{c}$, with CDF representations
\begin{align*}
F(z)=\sum_{j=1}^N r_j\one\{j\leq z\} \qquad &\forall z\in \RRR, \\
G(z)=\sum_{j=1}^N p_j\one\{j\leq z\} \qquad &\forall z\in \RRR.
\end{align*}
Note that $F(k)=R_k$ and $G^{-1}(P_k)=k$ when $k\in [N]$, and that $\RG=\{P_1,\ldots,P_{N-1}\}$.

The scoring rule $\GDS(\twocdots;G,\mu):\mathcal{P}_\mathrm{c}\times[N]\to[0,\infty)$, where $\mu=\sum_{k=1}^{N-1} w_k\delta_{P_k}$, is equitable with respect to $G$ and proper relative to $\mathcal{P}_\mathrm{c}$ by Theorems~\ref{thm:GDS is equitable} and \ref{thm:proper not strict}. Moreover, by \eqref{eq:gds} and Lemma~\ref{lem:F on integers},
\begin{align*}
\GDS(F,i;G,\mu)
&=\sum_{k=1}^{N-1} w_k(\one\{i<F^{-1}(P_k)\} - P_k)(\one\{G^{-1}(P_k)<F^{-1}(P_k)\} - \one\{G^{-1}(P_k)<i\}) \\
&=\sum_{k=1}^{N-1} w_k(\one\{F(i)<P_k\} - P_k)(\one\{F(k)<P_k\} - \one\{k<i\}) \\
&=S((r_1,\ldots,r_N), i),
\end{align*}
where $S$ is the scoring rule defined by (\ref{eq:proper equitable for cats}).
\end{proof}

\begin{lem}\label{lem:D(x,x')}
Suppose that $F,G\in\mathcal{P}$, $x'<x$ and that the expectations $\EEE_F \GDSd(x',Y;G,\mu)$ and $\EEE_F \GDSd(x,Y;G,\mu)$ are finite. Write $D(x,x')=\EEE_F \GDSd(x,Y;G,\mu)-\EEE_F \GDSd(x',Y;G,\mu)$.
\begin{enumerate}
\item[(i)] If $F(z)\leq G(z)$ for all $z$ in $[x',x)$ then $D(x,x')\leq 0$.
\item[(ii)]  If $F(z)\geq G(z)$ for all $z$ in $[x',x)$ then $D(x,x')\geq 0$.
\end{enumerate}
\end{lem}

\begin{proof}
From \eqref{eq:deterministic GDS},
\begin{align*}
&\EEE_F \GDSd(x,Y;G,\mu) \\
&\quad= \int_{\RG}\int_\RRR  (\one\{y<x\} - \alpha)(\one\{G^{-1}(\alpha) < x\} - \one\{G^{-1}(\alpha) < y\}) \,\dd F(y)\,\dd\mu(\alpha) \\
&\quad= \int_{\RG} \Big(\int_\RRR \one\{y<x\}\,\dd F(y) - \alpha \Big)\,\one\{G^{-1}(\alpha) < x\} \,\dd\mu(\alpha) \\
&\quad \qquad - \int_{\RG} \int_\RRR \one\{y<x\} \one\{G^{-1}(\alpha) < y\} \, \dd F(y) \,\dd\mu(\alpha) \\
&\quad \qquad + \int_{\RG} \alpha \int_\RRR \one\{G^{-1}(\alpha) < y\} \, \dd F(y) \,\dd\mu(\alpha) \\
&\quad= \int_{\RG} (F(x^-)-\alpha)\one\{G^{-1}(\alpha)<x\}\,\dd\mu(\alpha) \\
&\quad \qquad - \int_{\RG} \one\{G^{-1}(\alpha) < x\} (F(x^-)-F(G^{-1}(\alpha))) \,\dd\mu(\alpha) \\
&\quad \qquad + \int_{\RG} \alpha (1 - F(G^{-1}(\alpha))) \, \dd\mu(\alpha) \\
&\quad= \int_{\RG} \big( F(G^{-1}(\alpha)) - \alpha) \one\{G^{-1}(\alpha)<x\} + \alpha - \alpha F(G^{-1}(\alpha)) \big) \,\dd\mu(\alpha) \\
&\quad= \int_{\RG} \big( F(G^{-1}(\alpha))(1-\alpha)\one\{G^{-1}(\alpha)<x\} + \alpha(1-F(G^{-1}(\alpha))\one\{G^{-1}(\alpha)\geq x\} \big) \,\dd\mu(\alpha),
\end{align*}
where the interchange of integration in the first equality is justified by Fubini's theorem. Consequently,
\begin{align}
D(x,x')
&= \int_{\RG} \Big( F(G^{-1}(\alpha))(1-\alpha)(\one\{G^{-1}(\alpha) < x\} - \one\{G^{-1}(\alpha) < x'\}) \notag\\ 
&\qquad\phantom{\int_0^1 \Big(} + \alpha(1-F(G^{-1}(\alpha)))(\one\{G^{-1}(\alpha) \geq x\}- \one\{G^{-1}(\alpha) \geq x'\}) \Big) \,\dd\mu(\alpha) \notag\\
&= \int_{\RG} \Big( F(G^{-1}(\alpha))(1-\alpha) - \alpha(1-F(G^{-1}(\alpha)))\Big)\one\{ x'\leq G^{-1}(\alpha) < x\} \,\dd\mu(\alpha) \notag\\
&= \int_{\RG} (F(G^{-1}(\alpha)) - \alpha) \one\{x' \leq G^{-1}(\alpha) < x\} \,\dd\mu(\alpha). \label{eq:D(x,x')}
\end{align}

To prove (i), suppose that $F(z)\leq G(z)$ for all $z$ in $[x',x)$. Then
\begin{equation}\label{eq:D lemma}
F(G^{-1}(\alpha)) \leq G(G^{-1}(\alpha))=\alpha
\end{equation}
whenever $\alpha\in\RG$ and $x' \leq G^{-1}(\alpha) < x$, where the final equality in \eqref{eq:D lemma} is justified by Lemma~\ref{lem:alpha in RG}(i).
It follows from \eqref{eq:D(x,x')} that $D(x,x')\leq 0$.

To prove (ii), suppose that $F(z)\geq G(z)$ for all $z$ in $[x',x)$. Then $F(G^{-1}(\alpha))-\alpha \geq G(G^{-1}(\alpha)) - \alpha \geq 0$ whenever $\alpha\in\RG$ and $x' \leq G^{-1}(\alpha) < x$. It follows from \eqref{eq:D(x,x')} that $D(x,x')\geq 0$.
\end{proof}

\begin{proof}[Proof of Theorem~\ref{thm:optimal}]
The theorem follows from Definition~\ref{def:crossing point} and Lemma~\ref{lem:D(x,x')}.
\end{proof}

\begin{proof}[Proof of Proposition~\ref{prop:optimal fcst}]
Part (i) follows from the fact that in each case, $S$ is bounded. To prove part (ii), let $a_0$ and $b_0$ be the constants from Definition~\ref{def:global sharp}. We start with the decomposition of (\ref{eq:S2 decomp}) but with $\EEE_G$ replaced with $\EEE_F$. The first term on the right-hand side is finite. For the second term,
\begin{align*}
&\EEE_F[-\one\{x < Y\}\ln(1-G(Y))] \\
&\quad= \EEE_F[-\one\{\max(x,b_0) < Y\}\ln(1-G(Y))] + \EEE_F[-\one\{x < Y \leq \max(x,b_0)\}\ln(1-G(Y))] \\
&\quad\leq \EEE_G[-\one\{\max(x,b_0) < Y\}\ln(1-G(Y))] - \ln(1-G(\max(x,b_0))) \\
&\quad\leq \EEE_G[-\ln(1-G(Y))] - \ln(1-G(\max(x,b_0))) \\
&\quad\leq 1- \ln(1-G(\max(x,b_0))),
\end{align*}
where the first inequality is justified by the fact that $F(y)\geq G(y)$ and $y\mapsto -\ln(1-G(y))$ is increasing when $y > \max(x,b_0)$, and the third inequality is justified by Lemma~\ref{lem:finite expectations}. For the third term,
\begin{align*}
&\EEE_F[-\one\{x \geq Y\}\ln(G(Y))] \\
&\quad= \EEE_F[-\one\{\min(x,a_0) \geq Y\}\ln(G(Y))] + \EEE_F[-\one\{x \geq Y > \min(x,a_0)\}\ln(G(Y))] \\
&\quad\leq \EEE_G[-\one\{\min(x,a_0) \geq Y\}\ln(G(Y))] - \ln(G(\min(x,a_0))) \\
&\quad\leq \EEE_G[-\ln(G(Y))] -  \ln(G(\min(x,a_0))) \\
&\quad\leq 1- \ln(G(\min(x,a_0))),
\end{align*}
where the first inequality is justified by the fact that $F(y)\leq G(y)$ and $y\mapsto -\ln(G(y))$ is decreasing when $y\leq \min(x,a_0)$, and the third inequality is justified by Lemma~\ref{lem:finite expectations}. This completes the proof of part (ii).
\end{proof}

\begin{proof}[Proof of Proposition~\ref{prop:optimal seeps}]
For $i\in\{1,2,3\}$, let $e_i$ denote the expected SEEPS score given that category $C_i$ is forecast and the categorical outcome $Y_\mathrm{c}$ is distributed according to the predictive distribution $r=(r_1,r_2,r_3)$. That is, $e_i=\EEE_r[\SEEPS(C_i,Y_\mathrm{c})]$. 
Starting with the scoring matrix (\ref{eq:SEEPS}), we have
\begin{align*}
2e_1 &= \frac{r_2}{1 - p_1} + \left( \frac{1}{p_3} + \frac{1}{1 - p_1} \right) r_3, \\
2e_2 &= \frac{r_1}{p_1} + \frac{r_3}{p_3}, \\
2e_3 &= \left( \frac{1}{p_1} + \frac{1}{1 - p_3} \right) r_1 + \frac{r_2}{1 - p_3}. 
\end{align*}
The probability $r_2$ is eliminated from these equations using the fact that $r_1 + r_2 + r_3 = 1$. The first characterization in Proposition~\ref{prop:optimal seeps} is then obtained by rearranging and simplifying the inequalities $e_1\leq e_2$ and $e_1\leq e_3$. For example, it is straightforward to show that $e_1\leq e_2 \iff r_1\geq p_1$. The remaining characterizations are obtained similarly.
\end{proof}

\begin{lem}\label{lem:crossing point for normal}
Suppose that $G=\cN(m_1,\sigma_1^2)$ and $F=\cN(m_2,\sigma_2^2)$, where $0<\sigma_2<\sigma_1$. Then $F$ has exactly one crossing point $\chi$ on $\RRR$ with respect to $G$ and $\chi=(m_2\sigma_1-m_1\sigma_2)/(\sigma_1-\sigma_2)$.
\end{lem}

\begin{proof}
Since $0<\sigma_2<\sigma_1$ there must be at least one crossing point by the intermediate value theorem. Let $\Phi$ denote the CDF of $\cN(0,1)$. To find where the graphs of $F$ and $G$ intersect, we solve $\Phi((x-m_1)/\sigma_1)=\Phi((x-m_2)/\sigma_2)$ for $x$. Applying $\Phi^{-1}$ to both sides and rearranging for $x$ gives the unique solution $x=(m_2\sigma_1-m_1\sigma_2)/(\sigma_1-\sigma_2)$.
\end{proof}

\bibliography{equitability.bib}
\bibliographystyle{apalike-ejor}

\end{document}